\magnification=\magstep1
\input amstex
\documentstyle{amsppt}
\catcode`\@=11
\loadmathfont{rsfs}
\def\mycal{\mathfont@\rsfs}
\csname rsfs \endcsname
\catcode`\@=\active

\vsize=6.5in
\topmatter
\title STRONG RIGIDITY OF II$_1$ FACTORS
ARISING FROM MALLEABLE ACTIONS OF $w$-RIGID GROUPS, II
\endtitle
\author SORIN POPA \endauthor

\rightheadtext{Strong rigidity of factors, II}

\affil University of California, Los Angeles\endaffil

\address Math.Dept., UCLA, LA, CA 90095-155505\endaddress
\email popa\@math.ucla.edu\endemail

\thanks Supported in part by a NSF Grant 0100883.\endthanks

\abstract We prove that any isomorphism $\theta:M_0\simeq M$ of
group measure space II$_1$ factors, $M_0=L^\infty(X_0, \mu_0)
\rtimes_{\sigma_0} G_0$, $M=L^\infty(X, \mu) \rtimes_{\sigma} G$,
with $G_0$ an ICC group containing an infinite normal subgroup
with the relative property (T) of Kazhdan-Margulis (i.e. $G_0$
{\it w-rigid}) and $\sigma$ a Bernoulli shift action of some group
$G$, essentially comes from an isomorphism of probability spaces
which conjugates the actions with respect to some identification
$G_0 \simeq G$. Moreover, any isomorphism $\theta$ of $M_0$ onto a
``corner'' $pMp$ of $M$, for $p\in M$ an idempotent, forces $p=1$.
In particular, all group measure space factors associated with
Bernoulli shift actions of w-rigid ICC groups have trivial
fundamental group and any isomorphism of such factors comes from
an isomorphism    of the corresponding groups. This settles a
``group measure space version'' of Connes rigidity conjecture,
shown in fact to hold true in a greater generality than just for
ICC property (T) groups. We apply these results to ergodic theory,
establishing new strong rigidity and superrigidity results for
orbit equivalence relations.

\endabstract
\endtopmatter

\document
\heading 0. Introduction. \endheading

We continue in this paper the study of rigidity properties of
isomorphisms $\theta: M_0 \simeq M$ of crossed product II$_1$
factors initiated in ([Po4]), concentrating here on the ``group
measure space'' case, when the factors $M_0, M$ involved  come
from free ergodic measure preserving (m.p.) actions of groups on
probability spaces. Similarly to ([Po4]), we assume the ``source''
factor $M_0$ comes from an arbitrary free ergodic measure
preserving action $\sigma_0$ of a ``mildly rigid'' group $G_0$
(i.e. having a ``large'' subgroup with the relative property (T)
of Kazhdan-Margulis), while the ``target'' factor $M$ comes from
an action $\sigma$ satisfying good ``deformation+mixing''
properties (typically a Bernoulli shift). Our main result shows
that any isomorphism $\theta$  between such group measure space
factors $M_0, M$ essentially comes from an isomorphism of
probability spaces that conjugates the actions $\sigma_0, \sigma$,
with respect to some isomorphism of the groups $G_0, G$. In
particular, any isomorphism of the probability spaces that takes
the orbits of $\sigma_0$ onto the orbits of $\sigma$ (almost
everywhere), must come from a conjugacy of the actions $\sigma_0,
\sigma$. This establishes new strong rigidity and superrigidity
results in orbit equivalence ergodic theory, which are thus
obtained through a von Neumann algebra approach.

The factors $M_0, M$ come from the following data: $\sigma_0: G_0
\rightarrow \text{\rm Aut}(X_0,\mu_0)$, $\sigma: G \rightarrow
\text{\rm Aut}(X,\mu)$ are free ergodic measure preserving (m.p)
actions of countable discrete groups $G_0, G$ on standard
probability spaces $(X_0, \mu_0), (X,\mu)$;
$M_0=L^\infty(X_0,\mu_0) \rtimes_{\sigma_0} G_0$, $M =
L^\infty(X,\mu) \rtimes_{\sigma} G$ are their associated group
measure space (or crossed product) $\text{\rm II}_1$ factors, as
defined by Murray and von Neumann ([MvN1]). In addition to being
free and ergodic, $(\sigma_0, G_0), (\sigma,G)$ are assumed to
satisfy:

\vskip .1in \noindent $(*)$. $G_0$ is {\it w-rigid}, i.e. it
contains an infinite normal subgroup $H_0$ with the relative
property (T) of Kazhdan-Margulis (in other words,  $(G_0,H_0)$ is
a {\it property} (T) {\it pair}, see [Ma], [dHV]). We also
consider the class $w\Cal T_0$ of groups $G_0$ containing
subgroups  $H\subset G_0$ with the properties: $(a)$. $(G_0,H)$ is
a property (T) pair; $(b)$. $H$ is not virtually abelian; $(c)$.
$H$ satisfies a ``very weak'' normality condition in $G_0$ (see
$7.0.3$). The w-rigid groups will always be required ICC ({\it
infinite conjugacy class}), while the groups in the class $w\Cal
T_0$ need not be ICC.

Infinite property (T) groups are both w-rigid and in the class
$w\Cal T_0$. The groups $\Bbb Z^2 \rtimes \Gamma$, for $\Gamma
\subset SL(2, \Bbb Z)$ non-amenable (cf [Ka], [Ma], [Bu]) and the
groups $\Bbb Z^N \rtimes \Gamma$ for suitable actions of
arithmetic lattices $\Gamma \subset SU(n, 1), SO(n,1), n \geq 2$
(cf [Va], [Fe]) are all w-rigid, but not in the class $w\Cal T_0$.
Any product of a w-rigid (resp. $w\Cal T_0$) group with an
arbitrary group is still w-rigid (resp. $w\Cal T_0$). \vskip .1in
\noindent $(**)$. $\sigma$ is a Bernoulli shift action of an ICC
group $G$. Thus, $\sigma$ acts on the probability space $(X,
\mu)=\Pi_g (Y_0, \nu_0)_g$ by $\sigma_g((x_h)_h)=(x'_h)_h$, where
$x'_h=x_{g^{-1}h}, \forall h$. The key attributes of Bernoulli
shifts that we actually need are: a deformation property called
{\it sub s-malleability} (Definition 4.2 in [Po4]); a strong
mixing property called {\it clustering} (Definition 1.1 in this
paper). Thus, all statements below that require an action to be
Bernoulli, equally work with the assumption ``sub s-malleable and
clustering''.

\vskip .1in

The results we prove concern the structure of isomorphisms
$\theta:M_0 \simeq M$ between factors as above. More generally, we
consider isomorphisms between amplifications of $M_0, M$ by
positive real numbers ([MvN2]): If $N$ is a II$_1$ factor and
$t>0$ then the {\it amplification} of $N$ by $t$ is the
(isomorphism class of the) II$_1$ factor $N^t\overset\text{\rm
def} \to = pM_{n \times n}(N)p$, where $n \geq t$ and $p$ is a
projection in $M_{n \times n}(N)$ of normalized trace $\tau(p)$
equal to $t/n$. Thus, if $1\geq t
>0$ then $N^t$ is a ``corner'' $pNp$ of $N$, for $p\in N$ a
projection of trace $t$.

Any  {\it orbit equivalence} (OE) of $\sigma_0, \sigma$ (i.e. an
isomorphism of probability spaces $(X_0, \nu_0) \simeq (X,\mu)$
which takes onto each other the equivalence relations $\Cal
R_{\sigma_0}, \Cal R_{\sigma}$ given by the orbits of
$(\sigma_0,G_0), (\sigma,G)$) implements an isomorphism of $M_0,
M$ (cf. [Dy], [FM]). More generally, if $Y \subset X$ has measure
$\mu(Y)=t>0$ then any isomorphism $(X_0, \mu_0)\simeq (Y, \mu_Y)$
which takes $\Cal R_{\sigma_0}$ onto the equivalence relation
$\Cal R_\sigma^Y$ obtained by intersecting the orbits of $\sigma$
with $Y$ implements an isomorphism $M_0 \simeq M^t$. In general,
not all isomorphisms of group measure space factors come from OE
([CJ1]), and OE need not imply isomorphism of groups ([MvN2],
[Dy], [OW], [CFW]), much less conjugacy of actions. However, we
have:

\proclaim{0.1. Theorem (vNE Strong Rigidity)} Let $M_0,M$ be group
measure space $\text{\rm II}_1$ factors arising from actions
$(\sigma_0,G_0), (\sigma, G)$ satisfying $(*)$ resp. $(**)$. If
$\theta: M_0 \simeq M^t$ is an isomorphism of von Neumann
algebras, for some $1\geq t>0$, then $t=1$ and $\theta$ is of the
form $\theta = {\text{\rm Ad}}(u) \circ \theta^\gamma \circ
\theta^{\delta,\Delta}$, where: $u$ is a unitary element in $M$;
$\theta^\gamma \in {\text{\rm Aut}}(M)$ is implemented by a
character $\gamma$ of $G$; $\theta^{\delta,\Delta} :M_0 \simeq M$
is implemented by an isomorphism of groups $\delta : G_0 \simeq G$
and an isomorphism of probability spaces $\Delta : (X_0, \mu_0)
\simeq (X,\mu)$ satisfying $\sigma(\delta(h)) \circ \Delta =
\Delta \circ \sigma_0(h)$, $\forall h\in G_0$.
\endproclaim

This extremely rigid situation is in sharp contrast with the
amenable case, where by Connes' Theorem ([C3]) all group measure
space II$_1$ factors $L^\infty(X,\mu) \rtimes_\sigma G$ are
isomorphic, and by results of Dye ([Dy]), Ornstein-Weiss ([OW]),
Connes-Feldman-Weiss ([CFW]) all $(\sigma, G)$ are orbit
equivalent. Instead, Theorem 0.1 is in line with the rigidity
results established over the last 25 years in von Neumann algebra
theory ([C1,2], [CJ2], [Po8], [CoHa], [GoNe], [Po1-4]) and orbit
equivalence ergodic theory ([Zi], [GeGo], [Po8], [Ga1,2], [Fu1-3],
[MoSh], [Hj], [PoSa], [Po9]).

In fact, the terminology ``strong rigidity''  is borowed from
ergodic theory, where ``strong orbit rigidity'' ([Zi], [Fu1]), or
``OE strong rigidity'' ([MoSh]), designates results showing that
within a certain class of group actions orbital equivalence
automatically entails isomorphism of the groups and conjugacy of
the actions. While deriving the same type of conclusion, Theorem
0.1 assumes an even weaker equivalence of actions than orbit
equivalence (OE), namely {\it von Neumann equivalence} (vNE), i.e.
isomorphism of the associated group measure space von Neumann
algebras. Since Connes and Jones constructed examples of non-OE
actions that give rise to isomorphic factors ([CJ1]), vNE is
indeed (strictly) weaker than OE. Thus, Theorem 0.1 can be viewed
both as a ``strong rigidity''-type result in von Neumann algebra
theory, and as a new, stronger rigidity statement for ergodic
theory.

By applying this result to the case both actions are Bernoulli
shifts and both groups are ICC and either w-rigid or in $w\Cal
T_0$, one obtains a large class of group measure factors that can
be distinguished by the isomorphism class of their group-data:

\proclaim{0.2. Corollary} Let $G_i$ be ICC and either w-rigid or
in the class $w\Cal T_0$ and $\sigma_i : G_i \rightarrow
{\text{\rm Aut}} (X_i, \mu_i)$ a Bernoulli $G_i$-action, $i=0,1$.
Denote $M_i=L^\infty(X_i,\mu_i) \rtimes_{\sigma_i} G_i$ the
associated group measure space $\text{\rm II}_1$ factors and
$M_i^\infty=M_i \overline{\otimes} \Cal B(\ell^2\Bbb N)$ the
corresponding $\text{\rm II}_\infty$ factors, $i=0,1$. Then
$M_0^\infty \simeq M_1^\infty$ $\Leftrightarrow$ $M_0\simeq M_1$
$\Leftrightarrow$ $\sigma_0$ conjugate to
$\sigma_1$ with respect to some identification of
$G_0\simeq G$. Even more so, if $\theta: M^\infty_0 \simeq M^\infty_1$,
then $\theta$ is the amplification of an isomorphism $\theta_0 :
M_0 \simeq M_1$ of the form $\theta_0 = \text{\rm Ad}(u) \circ
\theta^\gamma \circ \theta^{\delta,\Delta}$, for some $u \in \Cal
U(M_1)$, $\gamma\in \text{\rm Char}(G_1)$ and $\delta:G_0 \simeq
G_1$, $\Delta:(X_0,\mu_0) \simeq (X_1,\mu_1)$ satisfying
$\sigma_1(\delta(h)) \Delta = \Delta \sigma_0(h)$, $\forall h\in
G_0$.
\endproclaim

At the Kingston AMS Summer School in 1980 A. Connes posed the
problem of showing that factors arising from property (T) ICC
groups are isomorphic iff the groups are isomorphic (see [C2]). He
formulated the question related to his discovery that property (T)
ICC groups $G$ give rise to group factors $L(G)$ with rigid
symmetry structure ([C1]). By now, several results in von Neumann
algebras ([CoHa], [Po8], [CSh]) and ergodic theory ([Zi], [CoZi],
[GeGo], [GoNa], [Fu1,2]) provide supporting evidence towards a
positive answer to the conjecture (see V.F.R. Jones' comments on
the ``higher rank lattice'' version of this problem in [J2]).

Since $L^\infty(X_i,\mu_i) \rtimes_{\sigma_i} G_i$, for $\sigma_i$
Bernoulli, can be viewed as a canonical ``group measure
space''-version of the group factor $L(G_i)$, Corollary 2 can be
regarded as an affirmative answer to a natural ``relative
variant'' of the conjecture. This is best emphasized by
re-formulating Corollary 0.2 as a ``group factor'' statement, by
using wreath product groups.  Thus, if we denote $\tilde{G}=\Bbb
Z^G \rtimes G$ (the wreath product with ``base'' $\Bbb Z$) then
the group factor $L(\tilde{G})$ is naturally isomorphic to the
group measure space factor corresponding to the Bernoulli shift
action $\sigma$ of $G$ on $\Bbb T^G$, and so we get:

\proclaim{0.3. Corollary} Let $G_i$ be w-rigid ICC groups and
denote $\tilde{G}_i=\Bbb Z^{G_i} \rtimes G_i$ the corresponding
wreath product, $i=0,1$. Then $L(\tilde{G_0})\simeq
L(\tilde{G}_1)$ iff $\tilde{G}_0 \simeq \tilde{G}_1$ and iff $G_0
\simeq G_1$.
\endproclaim

Corollary 0.2 shows that all group measure space factors
$L^\infty(X,\mu) \rtimes_{\sigma} G$ for $G$ w-rigid ICC and
$\sigma$ Bernoulli $G$-action have trivial fundamental group. It
also reduces the calculation of the outer automorphism group of
such a factor to the calculation of the commutant of $\sigma$ in
the set of automorphisms of the probability space $(X,\mu)$,
$\text{\rm Aut}_0(\sigma,G)\overset \text{\rm def}
\to=\sigma(G)'\cap \text{\rm Aut}(X,\mu)$:

\proclaim{0.4. Corollary} Let $G$ be ICC and either w-rigid or in
$w\Cal T_0$. If $\sigma: G \rightarrow \text{\rm Aut}(X,\mu)$ is a
Bernoulli $G$-action and $M=L^\infty(X,\mu) \rtimes_\sigma G$,
$M^\infty = M \overline{\otimes} \Cal B(\ell^2\Bbb N)$ then $M$
has trivial fundamental group, $\mycal F(M)=\{1\}$, and $\text{\rm
Out}(M^\infty)=\text{\rm Out}(M)= \text{\rm Aut}_0(\sigma,G)
\times (\text{\rm Char}(G) \rtimes \text{\rm Out}(G))$. Similarly,
$\mycal F(\Cal R_\sigma) = \{1\}$, $\text{\rm Out}(\Cal R_\sigma)
=\text{\rm Aut}_0(\sigma,G) \times \text{\rm Out}(G)$.
\endproclaim

The first examples of II$_1$ factors with trivial fundamental
group were obtained in ([Po2]; see also [Po3]). The key result
which allowed this computation was a ``strong rigidity''-type
result similar to Theorem 0.1, showing that any stable vNE of two
actions $(\sigma_0,G_0), (\sigma,G)$, with $\sigma_0$ satisfying a
relative property (T) in the spirit of Kazhdan-Margulis and
Connes-Jones ([CJ2]), and with $G$ satisfying Haagerup's compact
approximation property, comes from an OE of the actions. Thus, for
the class $\Cal H\Cal T$ of factors arising from actions
satisfying both properties, all OE invariants (such as Gaboriau's
cost or $\ell^2$-Betti numbers [Ga1,2]) are in fact isomorphism
invariants of the factors. Using the calculations of
$\ell^2$-Betti numbers in ([Ga2]) one derives that $L^\infty(\Bbb
T^2, \lambda) \rtimes \Bbb F_n = L(\Bbb Z^2 \rtimes \Bbb F_n)$
have trivial fundamental group and are mutually non-isomorphic for
$n = 2, 3, ...$. Thus, the approach for calculating fundamental
groups in ([Po2]) makes crucial use of Gaboriau's work in OE ergodic
theory.

In turn, the approach in this paper leads to a direct calculation
of fundamental group of factors, through purely von Neumann
algebra techniques, without using OE ergodic theory. Even more so,
it provides new calculations of fundamental groups of equivalence
relations, not covered by ([Ga1,2], [Fu1,2], [MoSh]).

We already pointed out that Theorem 0.1 implicitly gives an OE
strong rigidity result. We mention two more OE rigidity results
that we can derive from 0.1 (see Section 7 for more applications).
The first one concerns embeddings of equivalence relations and it
can be viewed as a Galois-type correspondence between w-rigid
subequivalence relations of a Bernoulli $G$-action  and w-rigid
subgroups of $G$:

\proclaim{0.5. Theorem (OE Strong Rigidity for Embeddings)} Let
$\sigma : G \rightarrow {\text{\rm Aut}} (X, \mu)$ be a  Bernoulli
shift action of an ICC group $G$. If $\sigma_0$ is a free ergodic
m.p. action of a w-rigid ICC group $G_0$ (or of a group $G_0\in
w\Cal T_0$) on $(X_0, \mu_0)$ and $\Delta_0 : (X_0, \mu_0) \simeq
(X,\mu)$ satisfies $\Delta_0(\Cal R_{\sigma_0})\subset \Cal
R_\sigma$ then there exists an isomorphism $\delta:G_0 \simeq
G_0'\subset G$ and $\alpha \in {\text{\rm Inn}}(\Cal R_\sigma)$
such that $\Delta=\alpha \circ \Delta_0$ satisfies $\Delta \circ
\sigma_0(h) = \sigma(\delta(h))\Delta$, $\forall h\in G_0$.
\endproclaim

If we restrict ourselves from w-rigid to Kazhdan groups, one can
deduce from Theorem 0.1 an ``OE rigidity'' result where all
conditions on the actions and on the groups involved are on just
``one side'', a type of result labelled ``OE superrigidity'' in
([MoSh]):

\proclaim{0.6. Theorem (OE Superrigidity)} Let $G$ be an ICC
property $(\text{\rm T})$ group and $\sigma: G \rightarrow
\text{\rm Aut}(X,\mu)$ a Bernoulli shift action. Let $G_0$ be any
group and $\sigma_0 : G_0 \rightarrow \text{\rm Aut}(X_0, \mu_0)$
any free m.p. action of $G_0$. If $\Delta_0: (X_0, \mu_0) \simeq
(Y, \mu_{|Y})$ satisfies $\Delta_0(\Cal R_{\sigma_0})=\Cal
R^Y_\sigma$, for some subset $Y\subset X$ of positive measure,
then $Y=X$ and there exist $\alpha \in \text{\rm Inn}(\Cal
R_\sigma)$, $\delta : G_0 \simeq G$, $\Delta : (X_0, \mu_0) \simeq
(X,\mu)$ satisfying $\sigma(\delta(h)) \Delta = \Delta
\sigma_0(h)$, $\forall h\in G_0$, such that $\Delta_0 = \alpha
\circ \Delta$. In particular, if $Y \subset X$ has measure $0\neq
\mu(Y) < 1$ then the equivalence relation $\Cal R^Y_\sigma$ cannot
be realized as orbits of a free action of a group.
\endproclaim

OE rigidity results were pioneered by Zimmer in the early 80's
(see [Zi]).  Furman revived this direction, establishing further,
sweeping OE rigidity results for classic actions of higher rank
lattices ([Fu1]). Most recently, Monod and Shalom proved surprising ``OE
strong rigidity and superrigidity'' results for a completely new class of
groups (e.g. products of torsion free non-elementary hyperbolic
groups) under very general ergodicity conditions on the actions
([MoSh]). The above applications bring additional insight to this
fascinating subject.

A few words on the proofs: It is the sub s-malleability (+ mixing)
condition on $\sigma$ and the ``weak property (T)'' (w-rigidity)
of $G_0$ that were used in ([Po4]) to prove that any
$\theta:M_0\simeq pMp$ as in Theorem 0.1 can be perturbed by an
inner automorphism of $pMp$ so that to take the group von Neumann
algebra $L(G_0)$, generated by $G_0$ in $M_0$, into a corner of
the group von Neumann algebra $L(G)$, generated by $G$ in $M$.
This preliminary rigidity result is crucial for the analysis in
this paper. It reduces the proof of Theorem 0.1 to settling the
case $M_0=pMp$, $p\in L(G)$, $L(G_0)\subset pL(G)p$. We in fact
prove a much more general statement, which only requires $\sigma$
clustering, letting $G, G_0$, $\sigma_0$ arbitrary:

\proclaim{0.7. Theorem (Criterion for Conjugacy of Actions)} Let
$\sigma : G \rightarrow {\text{\rm Aut}} (X, \mu)$, $\sigma_0 :
G_0 \rightarrow {\text{\rm Aut}}(X_0, \mu_0)$ be free, ergodic,
m.p. actions. Denote $A=L^\infty(X,\mu)$, $M=A \rtimes_\sigma G$,
$A_0 = L^\infty(X_0, \mu_0)$, $M_0 = A_0 \rtimes_{\sigma_0} G_0$.
Let $L(G)\subset M$ (respectively $L(G_0)\subset M_0$) be the von
Neumann subalgebra generated by the canonical unitaries
$\{u_g\}_{g\in G}\subset M$ implementing the action $\sigma$
(resp. $\{u^0_h\}_{h\in G_0}\subset M_0$ implementing $\sigma_0$).
Let also $p\in L(G)$ be a projection and assume:

$(a)$. $\sigma$ is clustering.

$(b)$. $pMp = M_0$,  with $L(G_0)$ contained in $pL(G)p$.

Then $p$ is central in $L(G)$, $\tau(p)^{-1}$ is an integer,
$pL(G)p = L(G_0)$ and there exist a normal subgroup $K \subset G$
with $|K|=\tau(p)^{-1}$, an ${\text{\rm Ad}}(G)$-invariant
character $\kappa$ of $K$, implementing a trivial $2$-cocycle on
$G/K$, and a unitary element $u\in pL(G)p$ such that:

\vskip .05in $(i)$. $p=|K|^{-1} \Sigma_{k \in K} \kappa(k) u_k$
and $u_kp=\kappa(k) p, \forall k\in K$.

\vskip .05in $(ii)$. $uA_0u^* = A^Kp$, where $A^K=\{a\in A \mid
\sigma_k(a)=a, \forall k\in K\}$.

\vskip .05in $(iii)$. $\{u_gp \mid g\in G\} = \{uu^0_hu^* \mid
h\in G_0\}$, modulo multiplication by scalars. More precisely,
there exist an isomorphism $\delta: G_0 \simeq G/K$, a lifting
$G/K \ni g \mapsto g' \in G$ and a map $\alpha: G_0 \rightarrow
\Bbb T$ such that $\text{\rm Ad}(u)(u^0_h)= \alpha(h)
u_{\delta(h)'} p$, $\forall h\in G_0$.

If in addition $G$ is assumed ICC then $K=\{e\}$, $p=1$ and
$\delta$ is an isomorphism $\delta:G_0\simeq G$.
\endproclaim

Note that this result holds true even in the hyperfinite case,
i.e. when $G_0, G$ are amenable. In fact, if we apply it to the
particular case $M_0=M=R$, $G=G_0$ abelian, $p=1$ and
$L(G)=L(G_0)$ then conclusion $(ii)$ alone gives an affirmative
answer to a recent conjecture of Neshveyev and St\o rmer in
([NeSt]), but for the more restrictive class of clustering actions
(e.g. Bernoulli shifts) rather than just weakly mixing, as
formulated in ([NeSt]).

This Conjugacy Criterion is stated as Theorem 5.1 in the text. Its
proof, which occupies most of the paper (sections 1-6), splits
into two parts, using rather different techniques: Part I consists
in proving a ``Cartan Conjugacy Criterion'' (Theorem 4.2 in the
text), showing that assumptions $(a), (b)$ above imply $A_0$ can
be unitary conjugate onto a corner of $A$. The proof, which takes
Sections 1-4, utilises (A.1 in [Po2]), ultrapower algebra
techniques and a careful ``asymptotic analysis'' of Fourier
expansions $\Sigma_g a_g u_g$ for elements in $A \rtimes_\sigma
G$. It is the only part that uses the {\it clustering} assumption
on $\sigma$ (see 1.1 below for the definition of clustering).

Part II of the proof of 0.7 consists in showing that if in
addition to $(a), (b)$ we also assume $A_0$ can be unitary
conjugate onto a corner of $A$, then $(i)-(iii)$ follow. This is
Theorem 5.2 in the text. In fact, instead of $(a)$ (i.e. the
clustering condition on $\sigma$), for this result we only assume
$\sigma$ and $\sigma_0$ mixing. The proof takes Sections 5 and 6
and uses ``local quantization'' techniques, in the spirit of (A.1
in [Po9]).

Our results leave unsolved the problem of showing that a group
measure space factor arising from a Bernoulli action of an ICC
property (T) group has ALL its Cartan subalgebras unitary
conjugate. Or at least showing that an isomorphism of group
measure space factors $M_0, M$ arising from actions
$(\sigma_0,G_0), (\sigma,G)$ is sufficient to imply that if $G$
has (T) then $G_0$ has (T). If true, then Theorem 0.1 would also
imply a ``vNE superrigidity''-type result for factors, similar to
the OE superrigidity in Theorem 0.6.

The present article is an outgrowth of a paper with the same title
and references (math.OA/0407137) that we have circulated since
July 2004 and in which we only proved that an isomorphism
$\theta:M_0 \simeq M^t$ with $M_0, M$ satisfying $(*), (**)$ comes
from an orbit equivalence of $\sigma_0, \sigma^t$. To derive from
that the triviality of the fundamental group and rigidity results
for factors, we had to make additional restrictions on the groups
and rely on OE rigidity results of Gaboriau ([Ga]), Monod-Shalom
([MoSh]) and ([PoSa]). All these applications are now consequences
of the ``vNE Strong Rigidity'' (0.1) and are shown is a greater
generality.

\heading 1. Clustering properties for actions
\endheading

The purpose of the first three sections is to lay down the
necessary technical background for the proof of the Cartan
Conjugacy Criteria in Section 4.

Thus, in this Section we define the ``clustering'' property for
actions of groups and a related notion of ``clustering
coefficients'' for sequences of elements in the corresponding
cross-product algebras. The considerations we make in this
Section, as well as in Sections 2 and part of 3, concern actions
of groups by trace preserving automorphisms on arbitrary finite
von Neumann algebras $(N, \tau)$ ($\tau$ denoting the fixed
faithful normal trace on $N$). Starting with 3.5 though, the
finite von Neumann algebra $(N, \tau)$ will always be abelian, a
fact that will be emphasized by using the generic notation $A$
instead of $N$ (possibly with indices).

\vskip .05in \noindent {\bf 1.1. Definition}. Let $\sigma: G
\rightarrow {\text{\rm Aut}} (N, \tau)$ be an action of a discrete
group $G$ on a finite von Neumann algebra $(N, \tau)$. A pair
$(N^0, \{B_n\}_n)$ consisting of a dense $^*$-subalgebra $N^0
\subset N$ and a decreasing sequence of von Neumann subalgebras
$\{B_n\}_n$ of $N$ is a {\it clustering resolution} for $\sigma$
if it satisfies the following conditions: \vskip .05in \noindent
$(1.1.1)$. For all $y \in N^0$ there exists $m$ such that
$E_{B_m}(y) = \tau(y) 1$; \vskip .05in \noindent $(1.1.2)$. For
all $m$ and all $y \in N^0$ there exists $F \subset G$ finite such
that $\sigma_g(y) \in B_m$, $\forall g \in G \setminus F$; \vskip
.05in \noindent $(1.1.3)$. For all $m \geq 1$ and all $g \in G$
there exists $n$ such that $\sigma_g(B_n) \subset B_m$.

\vskip .05in \noindent An action $\sigma$ is {\it clustering} if
it has a clustering resolution.

\proclaim{1.2. Proposition} $1^\circ$. If $\sigma$ is clustering
then for any $n \geq 1$ it satisfies the following $n$-mixing
condition:
$$
\underset g_1,..., g_n  \rightarrow \infty \to \lim |\tau (y_0
\sigma_{g_1}(y_1) ...
\sigma_{g_n}(y_n))-\tau(y_0)\tau(\sigma_{g_1}(y_1) ...
\sigma_{g_n}(y_n))| =0 \tag 1.2.1
$$
for all $y_0, y_1, ..., y_n \in N$.

$2^\circ$. If $\sigma$ is either a commutative or a Connes-St\o
rmer Bernoulli shift action then $\sigma$ is clustering.
\endproclaim
\vskip .05in \noindent {\it Proof}. $1^\circ$. By $(1.1.1)$, and
Kaplansky's density theorem, $\underset n \rightarrow \infty \to
\lim \|E_{B_n}(x)-\tau(x) 1\|_2 =0$, $\forall x\in N$. Thus, for
all $x \in N$ we have $\underset n \rightarrow \infty \to \lim
|\tau(E_{B_n}(x)x')-\tau(x)\tau(x')| =0$ uniformly in $x'\in
(N)_1$. Letting $x=y_0$, $x' = \sigma_{g_1}(y_1) ...
\sigma_{g_n}(y_n)$ and applying $(1.1.2)$, it follows that $$
\underset g_1,...,g_n \rightarrow \infty \to \lim
\tau(y_0\sigma_{g_1}(y_1) ... \sigma_{g_n}(y_n)) = \tau(y_0)
\tau(\sigma_{g_1}(y_1) ... \sigma_{g_n}(y_n)).
$$

$2^\circ$. Let $(\Cal N_0, \varphi_0)$ be a von Neumann algebra
with discrete decomposition. Assume $\sigma$ is the $(\Cal N_0,
\varphi_0)$-Bernoulli shift action of $G$ on $(\Cal N, \varphi)
={\overline{\underset g \in G \to \otimes}} (\Cal N_0,
\varphi_0)_g$ given by $\sigma_h(\otimes_g x_g)= \otimes_g x'_g$,
where $x'_g=x_{h^{-1}g}$. It also acts on the centralizer algebra
$N=\Cal N_\varphi$, by restriction.

Let $S_n$ be an increasing sequence of finite subsets of $G$ such
that $\cup_n S_n = G$ and denote by $\Cal B_n = {\underset g \in
S_n^c \to {\overline{\otimes}}} (\Cal N_0, \varphi_0)_g$, where
$S_n^c = G \setminus S_n$. We claim that $N^0= \Cal N_\varphi \cap
\otimes_g (\Cal N_0, \varphi_0)_g \subset N$ (algebraic tensor
product) and $B_n = \Cal N_\varphi \cap \Cal B_n, n\geq 1,$ give a
clustering resolution for the action $\sigma$ of $G$ on $N=\Cal
N_\varphi$. Indeed, condition $(1.1.1)$ is clearly satisfied by
the definition and by the fact that $gh \in S_n^c$ as $g
\rightarrow \infty$, $\forall h\in G$.

To check $(1.1.2)$ note that if $y\in \Cal N_\varphi \cap
{\overline{\underset g \in S_n \to \otimes}} (\Cal N_0,
\varphi_0)_g$, for some large enough $n$, then for any $m$,
$\sigma_g(y) \in B_m$ as $g \rightarrow \infty$.

To show that $(1.1.3)$ is verified,
consider the weakly dense subalgebra
$B_n^0=B_n \cap N^0$ of $B_n$ and note that since $\cup_n
S_n = G$, we have $g^{-1}S_k \subset S_n$ (equivalently
$g(S_n^c) \subset S_k^c$) as $n \rightarrow \infty$. This implies that
$\sigma_g(B_n^0) \subset B_k$ as $n \rightarrow \infty$. But for each
$n$ for which we have the above inclusion we also have
$\sigma_g(B_n) \subset B_k$. \hfill Q.E.D.

\vskip .05in
We denote by $M = N \rtimes_\sigma G$ the cross-product von
Neumann algebra associated to an action $\sigma : G \rightarrow
{\text{\rm Aut}}(N, \tau)$ and by $\{u_g\}_g \subset M$ the
canonical unitaries implementing the action $\sigma$. Also, we
denote by $L(G)$ the von Neumann subalgebra of $M$ generated by
the unitaries $\{u_g\}_g$.

Let $\omega$ be a free ultrafilter on $\Bbb N$, also viewed as a
point in the Stone-Cech compactification of $\Bbb N$, $\omega \in
\overline{\Bbb N} \setminus \Bbb N$. Note that a subset $V\subset
\Bbb N$ in the ultrafilter $\omega$ corresponds in this latter
case to a neighborhood of the point $\omega$ in the compact set
$\overline{\Bbb N}$.

As usual, given a finite von Neumann algebra $(B, \tau)$ we denote
by $(B^\omega, \tau)$ its {\it ultrapower} algebra, i.e. the
finite von Neumann algebra $\ell^\infty(\Bbb N, B)/\Cal I_\omega$,
where $\Cal I_\omega$ is the ideal of elements $x=(x_n)_n \in
\ell^\infty(\Bbb N, B)$ with $\tau(x^*x) = 0$, the trace
$\tau=\tau_\omega$ on $\ell^\infty(\Bbb N, B)/\Cal I_\omega$ being
defined by $\tau((b_n)_n) = \underset n \rightarrow \omega \to
\lim \tau(b_n)=0$. For basic properties and results on ultrapower
algebras see ([McD]).

We often identify $B$ with the subalgebra of constant sequences in
$B^\omega$.

\vskip .1in \noindent {\bf 1.3. Definition}. Let $\sigma : G
\rightarrow \text{\rm Aut}(N, \tau)$ be a properly outer action
with clustering resolution $(N^0, \{B_k\}_k)$. A sequence $(x_n)_n
\in \ell^\infty(\Bbb N, M)$ has {\it clustering coefficients} with
respect to $(N^0, \{B_n\}_n)$ (and to $\omega$) if $\forall
\varepsilon > 0$ and $\forall m\geq 1$, $\exists V \in \omega$
such that if $e_k$ denotes the orthogonal projection of $L^2(M)$
onto $L^2(\Sigma_g B_k u_g)$ then $\|e_{m}(\hat{x}_n) -
\hat{x}_n\|_2 \leq \varepsilon$, $\forall n \in V$. Note that if
$x=(x_n)_n$ has clustering coefficients then any element in $x +
\Cal I_\omega$ also has clustering coefficients. Thus, this
property passes to elements in $M^\omega = \ell^\infty(\Bbb N,
M)/I_\omega$. \vskip .1in \noindent {\bf 1.3.' Notations}. We
denote by $\Cal X$ the set of sequences $(x_n)_n\in
\ell^\infty(\Bbb N, M)$ with clustering coefficients. We still
denote by $\Cal X$ the image in $M^\omega$ of the set of sequences
with clustering coefficients. Also, we denote by $\Cal Y$ the set
of elements $(x_n)_n$ in $M^\omega$ with the property that
$\underset n \rightarrow \omega \to \lim\| E_N(x_nu_g)\|_2=0$,
$\forall g\in G$. \vskip .05in We denote by $N^\omega
\rtimes_\sigma G$ the von Neumann subalgebra of $M^\omega$
generated by $N^\omega$ and the canonical unitaries $u_g, g \in G$
(regarded as constant sequences in $M^\omega$). Also, we view the
group von Neumann algebra $L(G)=\{u_g\}_g''$ as a subalgebra of
$M^\omega$, via the inclusion $M \subset M^\omega$.

\proclaim{1.4. Lemma} $1^\circ$. $\Cal X, \Cal Y$ are vector
spaces.

$2^\circ$. If $\{x^k\}_k \subset \Cal X$ is convergent in the norm
$\|\quad \|_2$ to an element $x \in M^\omega$ then $x\in \Cal X$.
Also, the unit ball of $\Cal X$ (resp. $\Cal Y$) is complete in
the norm $\| \quad \|_2$ given by the trace $\tau$ on $M^\omega$.
\endproclaim
\vskip .05in \noindent {\it Proof}. $1^\circ$ is trivial by the
definitions.

Since $M^\omega$ is a II$_1$ factor ([McD]), $(M^\omega)_1$ is
complete in the norm $\| \quad \|_2$. Thus, the completeness of
$(\Cal X)_1$ follows from the first part of $2^\circ$.

To prove the first part of $2^\circ$, let $\varepsilon > 0$ and $m
\geq 1$. Then there exists $k$ such that $\|x-x^k\|_2 <
\varepsilon/3$. Thus, if $x^k=(x^k_n)_n$, $x=(x_n)_n$, then there
exists a neighborhood $V'$ of $\omega$ such that $\|x^k_n - x_n
\|_2 \leq \varepsilon/3$, $\forall n \in V'$. But since $x^k \in
\Cal X$, it follows that there exists a neighborhood $V''$ of
$\omega$ such that $\|e_m(\hat{x}^k_n)-\hat{x}^k_n\|_2 \leq
\varepsilon/3$, $\forall n \in V''$. Hence, if we denote $V=V'\cap
V''$ then for $n \in V$ we have:
$$
\|e_m(\hat{x}_n) - \hat{x}_n \|_2 \leq \|e_m(\hat{x}_n) -
e_m(\hat{x}^k_n)\|_2
$$
$$
+ \|e_m(\hat{x}^k_n) - \hat{x}^k_n\|_2 + \|\hat{x}^k_n -
\hat{x}_n\|_2 \leq 3 \varepsilon/3 = \varepsilon.
$$
This shows that $x \in \Cal X$.

The proof that $\Cal Y \cap (M^\omega)_1$ is complete is similar.
\hfill Q.E.D.

\proclaim{1.5. Lemma} $\Cal Y= M^\omega \ominus (N^\omega \rtimes
G)$. Equivalently, if $x\in M^\omega$ then $x\in \Cal Y$ iff
$E_{N^\omega \rtimes G}(x)=0$. Also, $L(G)^\omega \cap \Cal Y =
L(G)^\omega \ominus L(G)$.
\endproclaim
\vskip .05in \noindent {\it Proof}. This is trivial by the
definitions. \hfill Q.E.D.

\heading 2. Multiplicative properties of clustering sequences
\endheading

Unless otherwise specified, throughout this section $\sigma$
denotes a clustering action of $G$ on $(N, \tau)$ with clustering
resolution $(N^0, \{B_n\}_n)$, as in Definition 1.1. Also, $\Cal
X$ denotes its associated set of sequences with clustering
coefficients, as in 1.3'. We prove here some multiplicative
properties for $\Cal X$.

\proclaim{2.1. Lemma} $1^\circ$. $\Cal X \cap N^\omega$ is equal
to the von Neumann algebra $\cap_m B_m^\omega$. Also,
$E_{N^\omega}(\Cal X) \subset \Cal X$.

$2^\circ$. $\Cal X$ is a $L(G)-L(G)^\omega$ bimodule.
\endproclaim
\vskip .05in \noindent {\it Proof}. $1^\circ$. We clearly have
$\cap_m B_m^\omega \subset \Cal X \cap N^\omega$, by the
definitions. Conversely, if $x=(x_n)_n \in N^\omega$ belongs to
$\Cal X$ then for any $m \geq 1$ and any $\varepsilon > 0$ there
exists a neighborhood $V$ of $\omega$ such that $\|E_{B_m}(x_n) -
x_n\|_2 \leq \varepsilon$, $\forall n \in V$. But this implies
$\|E_{B_m^\omega}(x)-x\|_2 \leq \varepsilon$ and since
$\varepsilon > 0$ was arbitrary, this shows that $x\in
B_m^\omega$. Since $m$ was arbitrary as well, $x \in \cap_m
B_m^\omega$.

Since for $x=(x_n)_n \in M^\omega$ we have
$E_{N^\omega}(x)=(E_N(x_n))_n$, the last part is now clear by the
first part and the definitions.

$2^\circ$. Let $x=(x_n)_n \in \Cal X$ and $y=(y_n)_n \in
L(G)^\omega$, with $y_n \in (L(G))_1, \forall g$. By the definition
1.3, for any $m \geq 1$ and $\varepsilon > 0$ there exists a
neighborhood $V$ of $\omega$ and $\{\xi_n\}_n \subset L^2(\Sigma_g B_mu_g)$
such that $\|\hat{x}_n -\xi_n \|_2 \leq \varepsilon$, $\forall n
\in V$. But then $\eta_n = \xi_n y_n$ belongs to $L^2(\Sigma B_m
u_g)$ and $\|\hat{x_ny_n} - \eta_n \|_2 \leq \varepsilon$,
$\forall n \in V$. Thus, $xy \in \Cal X$.

To prove that $\Cal X$ is a left $L(G)$-module, by $1.4$ it is
sufficient to show that if $x \in \Cal X$ then $u_gx \in \Cal X$,
$\forall g \in G$. By Definition 1.3, for any $\varepsilon > 0$
and any $m$ there exist $m_0$, $V \in \omega$ and $\xi'_n \in
L^2(\Sigma_g B_{m_0} u_g)$ such that $\sigma_g(B_{m_0}) \subset
B_m$ and $\|\hat{x}_n-\xi'_n\|_2 \leq \varepsilon$, $\forall n\in
V$. But then we have $u_g \xi_n' \in L^2(\Sigma_g B_m u_g)$. Thus,
$\|u_g\hat{x}_n - u_g \xi_n'\|_2 \leq \varepsilon$, $\forall n \in
V$. This shows that $(u_gx_n)_n \in \Cal X$. \hfill Q.E.D.

\proclaim{2.2. Lemma} The unitaries $u_g, g\in G$, normalize $\Cal
X \cap N^\omega = \cap_m B_m^\omega$ and $\Cal X \cap (N^\omega
\rtimes G)$ is equal to the von Neumann algebra $(\cap_m
B_m^\omega) \rtimes G$, generated by $\cap_m B^\omega_m$ and
$\{u_g\}_g$. Also, $E_{N^\omega \rtimes G}(\Cal X) \subset \Cal
X$.
\endproclaim
\vskip .05in \noindent {\it Proof}. By $(1.1.3)$ and Definition
1.3, we clearly have $u_g (\cap_m B_m^\omega) u_g^* = \cap_m
B_m^\omega$. By 2.1.2$^\circ$, the $^*$-algebra generated by $\cap
B_m^\omega$ and $\{u_g\}_g$ is contained in $\Cal X$. Conversely,
if $x=(x_n)_n \in \Cal X \cap (N^\omega \rtimes G)$ then
$x=\Sigma_g E_{N^\omega}(xu_g^*)u_g$ and by 2.1.1$^\circ$ each
$E_{N^\omega}(xu_g^*), g\in G,$ follows in $\Cal X\cap N^\omega =
\cap_m B_m^\omega$. Thus, $x \in (\cap_m B_m^\omega) \rtimes G$.

To see that if $x\in \Cal X$ then $E_{N^\omega \rtimes G}(x)\in
\Cal X$, note that $E_{N^\omega}(xu_g^*)u_g \in (\cap_m
B_m^\omega) \rtimes G$, so by taking finite sums and limits we get
$E_{N^\omega \rtimes G}(x) \in L^2((\cap_m B_m^\omega) \rtimes
G)$, implying that $E_{N^\omega \rtimes G}(x) \in (\cap_m
B_m^\omega) \rtimes G \subset \Cal X $. \hfill Q.E.D.

\proclaim{2.3. Corollary} If $x \in \Cal X$, then $x-E_{N^\omega
\rtimes G}(x) \in \Cal X \cap \Cal Y.$
\endproclaim
\vskip .05in \noindent {\it Proof}. By the last part of 2.2 we
have $x-E_{N^\omega \rtimes G}(x) \in \Cal X$ while by 1.5 we have
$x-E_{N^\omega \rtimes G}(x) \in \Cal Y$.
\hfill Q.E.D.

\proclaim{2.4. Lemma} $\Cal X \cap \Cal Y$ is a right $M$-module,
i.e., if $x\in \Cal X \cap \Cal Y$ then $xM \subset \Cal X \cap
\Cal Y$. In particular, $(L(G)^\omega \ominus L(G))M \subset \Cal
X$.
\endproclaim
\vskip .05in \noindent {\it Proof}. By $1.4$ it is sufficient to
prove that $x y \in \Cal X\cap \Cal Y$ for $y \in N^0\cup
\{u_g\}_g$. The case $y = u_g$ is trivial by the definitions. Let
$y \in N^0, \|y\|\leq 1$. Then clearly $x y \in \Cal Y$. To show
that $x y \in \Cal X$ as well, let $m \geq 1$, $\varepsilon
> 0$ . By Definition 1.1 there exists $F \subset G$
finite such that $\sigma_g(y) \in B_m$, $\forall g \in G \setminus
F$. On the other hand, by 1.3, since $x \in \Cal X$ there exists a
neighborhood $V$ of $\omega$ so that $\|e_m
(\hat{x}_n)-\hat{x}_n\|_2 \leq \varepsilon/2$, $\forall n\in V$.
Also, since $x \in \Cal Y$, $V$ can be chosen so that if $x_n =
\Sigma_g y_{n,g} u_g$, with $y_{n,g} \in N$, then $\|y_{n,g}\|_2
\leq \varepsilon/3|F|^{1/2}$, $\forall n \in V$. All this entails
the following estimates:

$$
\|e_m(\hat{x_n y}) - \hat{x_n y} \|_2^2 = \Sigma_{g \in G}
\|E_{B_m}(y_{n,g} \sigma_g(y))- y_{n,g}\sigma_g(y)\|_2^2
$$
$$
=\Sigma_{g \in F} \|E_{B_m}(y_{n,g} \sigma_g(y))-
y_{n,g}\sigma_g(y)\|_2^2 $$ $$
+\Sigma_{g \in G \setminus F}
\|E_{B_m}(y_{n,g} \sigma_g(y))- y_{n,g}\sigma_g(y)\|_2^2
$$
$$
\leq 4 \Sigma_{g \in F} \|y_{n,g}\|_2^2 + \Sigma_{g \in G
\setminus F} \|(E_{B_m}(y_{n,g})- y_{n,g}) \sigma_g(y)\|_2^2
$$
$$
< \varepsilon^2/2 + \|e_m(\hat{x}_n) - \hat{x}_n\|_2^2 <
\varepsilon^2.
$$

Thus $xy \in \Cal X$.
\hfill Q.E.D.

\vskip .05in While by 2.4 $\Cal X\cap \Cal Y$ is invariant to
multiplication by $M$ from the right, the next result shows that
multiplication from the left by elements in $M \ominus L(G)$ takes
the space $\Cal X\cap \Cal Y$ perpendicular to $\Cal X$:

\proclaim{2.5. Lemma} $(M \ominus L(G)) (\Cal X \cap \Cal Y) \perp
\Cal X$.
\endproclaim
\vskip .05in \noindent {\it Proof}. Let $x=(x_n)_n \in \Cal X \cap
\Cal Y$, $y \in M$ with $E_{L(G)}(y)=0$. We have to show that
$\underset n \rightarrow \omega \to \lim \tau(z_n^* y x_n)=0$,
$\forall z=(z_n)_n \in \Cal X$.

By linearity and Kaplansky's density theorem it is clearly
sufficient to prove the statement for $y$ of the form $y = y_gu_g$
for some $g \in G$ and $y_g \in N^0, \tau (y_g)=0$. Also, since by
2.1 we have $u_g (\Cal X \cap \Cal Y) \subset \Cal X \cap \Cal Y$,
by replacing $(x_n)_n$ by $(u_gx_n)_n$ it follows that it is in
fact sufficient to prove the statement for $y \in N^0$ with
$\tau(y)=0$. In addition, we may assume $\|y\| \leq 1$, $\|x\|
\leq 1, \|z\| \leq 1$.

Since $e_m(\hat{y})=\hat{E_{B_m}(y)}$, we have $e_m(y\xi) =
E_{B_m}(y)\xi$ for any $\xi \in L^2(\Sigma_g B_mu_g)$. As $y \in
N^0$ and $\tau (y)=0$, by $(1.1.1)$ there exists $m$ such that
$E_{B_m}(y)=0$. Let $\varepsilon > 0$. Since $x, z \in \Cal X$,
there exists a neighborhood $V$ of $\omega$ so that
$\|e_m(\hat{x}_n) - \hat{x}_n\|_2 \leq \varepsilon/2$ and
$\|e_m(\hat{z}_n) - \hat{z}_n\|_2 \leq \varepsilon/2$, $\forall
n\in V$. But then, for each $n \in V$ we have the estimates:

$$
|\tau(z_n^*yx_n)| \leq |\langle y e_m(\hat{x}_n), e_m(\hat{z}_n)
\rangle|
$$
$$
+ \|\hat{z}_n - e_m(\hat{z}_n)\|_2 \|yx_n\|_2 +
\|e_m(\hat{z}_n)\|_2 \|y\hat{x}_n - ye_m(\hat{x}_n)\|_2
$$
$$
\leq |\langle y e_m(\hat{x}_n), e_m(\hat{z}_n) \rangle| +
\varepsilon = \varepsilon,
$$
where for the last equality we have used
$$
\langle y e_m(\hat{x}_n), e_m(\hat{z}_n) \rangle = \langle
e_m(ye_m(\hat{x}_n)), e_m(\hat{z}_n) \rangle = \langle E_{B_m}(y)
e_m(\hat{x}_n), e_m(\hat{z}_n) \rangle = 0.
$$
\hfill Q.E.D.

\vskip .05in

Since $L(G)^\omega \ominus L(G)\subset \Cal X \cap \Cal Y$, Lemmas
2.4, 2.5 imply that $(M\ominus L(G)) (L(G)^\omega \ominus L(G))
\perp (L(G)^\omega \ominus L(G))M$. This orthogonality holds in
fact true in a larger generality:

\proclaim{2.6. Lemma} Let $\sigma$ be a free, mixing action of a
discrete group $G$ on a finite von Neumann algebra $(N,\tau)$ and
denote $M=A \rtimes_\sigma G$. If $(x_n)_n, (y_n)_n \in L(G)$ are
bounded sequences that weakly converge to $0$, then $\underset n
\rightarrow \infty \to \lim \tau(x_n^*a^*y_nb)=0$ for all $a,b \in
M\ominus L(G)$. Thus, $(M\ominus L(G))(L(G)^\omega \ominus L(G))
\perp (L(G)^\omega \ominus L(G))M$.
\endproclaim
\vskip .05in \noindent {\it Proof}. We may assume $\|x_n\|,
\|y_n\| \leq 1, \forall n$. By linearity and Kaplansky's density
theorem, it is sufficient to prove the statement for $a= a_0u_g, b
= b_0u_h$, where $g,h \in G$ and $a_0, b_0 \in N$. Also, since
$(u_gx_n)_n, (y_nu_h)_n$ are still converging weakly to $0$, it
follows that it is in fact sufficient to prove the statement for
$a,b \in N$ with $\tau(a)=\tau(b)=0$, $\|a\|, \|b\| \leq 1$. For
such $a,b$, if we let $x_n = \Sigma_g c^n_g u_g$, $y_n = \Sigma_g
d^n_g u_g$, $c^n_g, d^n_g \in \Bbb C$, be the Fourier expansions,
then $\tau(x_n^*a^*y_n) = \Sigma_{g\in G} c^n_g\overline{d^n_g}
\tau(a^*\sigma_g(b))$, with the sum being $\ell^1(G)$ convergent.

Let $\varepsilon > 0$. Since $\sigma$ is mixing and $\tau(a)=0$,
there exists a finite subset $F \subset G$ such that
$|\tau(a^*\sigma_g(b))| \leq \varepsilon/2$, $\forall g\in G
\setminus F$. By the weak convergence of the sequences $(x_n)_n,
(y_n)_n$, there exists $n_0 \geq 1$ such that $\Sigma_{g\in F}
|c^n_g \overline{d^n_g}| \leq \varepsilon/2$, $\forall n \geq
n_0$. Altogether, for $n \geq n_0$ we get the estimates:
$$
|\tau(x_n^*a^*y_n)|=|\Sigma_{g\in G} c^n_g\overline{d^n_g}
\tau(a^*\sigma_g(b))|
$$
$$
\leq \Sigma_{g\in F} |c^n_g \overline{d^n_g}| \|a\| \|b\| +
\Sigma_{g \not\in F} |c^n_g \overline{d^n_g}| |\tau(a^*
\sigma_g(b))|
$$
$$
\leq \varepsilon/2 + (\varepsilon/2)\Sigma_{g \not\in F} |c^n_g
\overline{d^n_g}| \leq \varepsilon.
$$
\hfill Q.E.D.

\heading 3. Clustering algebras associated to pairs of actions
\endheading

Let $\sigma : G \rightarrow \text{\rm Aut}(N, \tau)$ be a properly
outer clustering action with clustering resolution $(N^0,
\{B_n\}_n)$, as in Sections 1, 2. In this Section we assume that
for a non-zero projection $p\in L(G)={\{u_g\}_g}''$ the reduced
algebra $M_0=pMp$ has another cross-product decomposition $M_0 =
N_0 \rtimes_{\sigma_0} G_0$, for some properly outer action
$\sigma_0 : G_0 \rightarrow {\text{\rm Aut}}(N_0, \tau)$.

\vskip .1in \noindent {\bf 3.1. Notation}. We denote $\Cal X_0=
\Cal X \cap (N_0'\cap N_0^\omega)$. Note that by $1.4.1^\circ$,
$\Cal X_0$ is a vector subspace of $pM^\omega p$. We show that in
fact :

\proclaim{3.2. Theorem } $\Cal X_0$ is a von Neumann subalgebra of
$p(N^\omega \rtimes G)p$.
\endproclaim
\vskip .05in \noindent {\it Proof}. We first prove that $\forall
\varepsilon > 0$, $\exists y_0 \in \Cal U(N_0)$ such that
$\|E_{L(G)}(y_0)\|_2 \leq \varepsilon$. Indeed, for if not then
$\exists \varepsilon_0 > 0$ such that if we denote by $(\langle M,
e_{L(G)} \rangle, Tr)$ the basic construction algebra
corresponding to the inclusion $L(G) \subset M$, then
$$
Tr(e_{L(G)}ye_{L(G)}y^*) = \|E_{L(G)}(y)\|^2_2 \geq
\varepsilon^2_0, \forall y \in \Cal U(N_0).
$$

Taking weak limits of convex combinations of elements of the form
$ye_{L(G)}y^*$, for $y \in \Cal U(N_0)$, it follows that
$Tr(e_{L(G)} a) \geq \varepsilon_0$ for any $a \in K=
\overline{\text{\rm co}}^w \{ye_{L(G)}y^* \mid y \in \Cal U(N_0)
\}$. In particular, $Tr(e_{L(G)}a_0) \geq \varepsilon_0^2$ for
$a_0$ the unique element of minimal norm $\| \cdot \|_{2,Tr}$ in
$K$. Thus, $a_0\neq 0$. But $ya_0 y^* \in K$ and
$\|ya_0y^*\|_{2,Tr} = \|a_0\|_{2,Tr}$, $\forall y\in \Cal U(N_0)$,
so by the uniqueness of $a_0$ we have $ya_0y^*=a_0, \forall y \in
\Cal U(N_0)$, showing that $a_0 \in N_0'\cap \langle M, e_{L(G)}
\rangle$, $0 \leq a_0 \leq 1, Tr(a_0) \leq Tr(e_{L(G)})=1$, $a_0
\neq 0$.

By (Theorem 2.1 in [Po4]), this implies there exists a non-zero
projection $q\in N_0$, a projection $p\in L(G)$, a unital
isomorphism $\psi$ of $qN_0q$ into $pL(G)p$ and a partial isometry
$v\in M$ such that $vv^* = q$, $v^*v \in \psi(qN_0q)'\cap pL(G)p$
and $x v = v\psi(x)$, $\forall x\in qN_0q$. Moreover, by shrinking
$q$ if necessary, we may assume the central trace of $q$ in $N_0$
is a scalar multiple of a central projection of $N_0$ (any
projection in $N_0$ majorizes a projection with this latter
property).

But $\sigma$ is mixing and $\psi(qN_0q)$ is a subalgebra without
atoms in $L(G)$, so by (Theorem 3.1 in [Po4]) it follows that
$\psi(qN_0q)'\cap pL(G)p \subset L(G)$. Thus, $p_0=v^*v \in L(G)$
and $v^*N_0v \subset p_0L(G)p_0$. Moreover, by applying again (3.1
in [Po4]) it follows that the normalizer of $v^*N_0v$ in $p_0Mp_0$
is included in $p_0L(G)p_0$. This is a contradiction, since by
(3.5.2$^\circ$ in [Po4]) the normalizer of $v^*N_0v$ in $p_0Mp_0$
generates all $p_0Mp_0$ (because $qN_0q$ is regular in $q M q$,
due to the fact that $N_0$ is regular in $M$).

Let now $x=(x_n)_n \in \Cal X_0\subset \Cal X$ and denote $x'=x -
E_{N^\omega \rtimes G}(x)$. By Lemma 2.3, $x'$ lies in $\Cal X
\cap \Cal Y$. Let $\varepsilon > 0$ and choose $y_0 \in \Cal
U(N_0)$ such that $\|E_{L(G)}(y_0)\|_2 \leq \varepsilon$. Since
$[y_0, x]=0$ and $y_0 \in M \subset N^\omega \rtimes G$, it
follows that $[y_0, E_{N^\omega \rtimes G}(x)]=0$. Thus, $[y_0,
x']=0$.

But by Lemma 2.4 we have $x'y_0 \in \Cal X$, while by Lemma 2.5 we
have $(y_0-E_{L(G)}(y_0)) x' \perp \Cal X$. In particular,
$(y_0-E_{L(G)}(y_0)) x'$ is perpendicular to $x'y_0$. Thus, the
commutation relation $y_0x'=x'y_0$ entails the equation
$$
x'y_0 - (y_0-E_{L(G)}(y_0)) x' = E_{L(G)}(y_0)x',
$$
with the vectors $x'y_0, (y_0-E_{L(G)}(y_0)) x'$ perpendicular. By
Pythagora's Theorem, and using that $y_0$ is a unitary element, we
get

$$
\|x'\|_2 = \|x'y_0\|_2 \leq \|E_{L(G)}(y_0)x'\|_2 \leq
\varepsilon.
$$
Since $\varepsilon > 0$ was arbitrary, $x'=0$. Thus $x =
E_{N^\omega \rtimes G}(x) \in N^\omega \rtimes G$. But this also
shows that $\Cal X_0 = (\Cal X \cap (N^\omega \rtimes G)) \cap
(N_0'\cap N_0^\omega)$. Since by Lemma 2.2 $\Cal X \cap (N^\omega
\rtimes G)$ is a weakly closed $^*$-algebra, this shows that $\Cal
X_0$ is a von Neumann algebra as well. \hfill Q.E.D.

\vskip .05in
The next Lemma and its Corollary will be needed in Section 4.
\proclaim{3.3. Lemma} Given any $\varepsilon > 0$ there exists a
finite subset $F \subset G$ such that if $p_F$ denotes the
orthogonal projection of $L^2(M^\omega)$ onto $L^2(\Sigma_{g\in F}
N^\omega u_g)$ then for any $x \in (\Cal X_0)_1$ we have
$\|\hat{x}-p_{F}(\hat{x})\|_2 \leq \varepsilon$.
\endproclaim
\vskip .05in \noindent {\it Proof}. Assume this is not the case.
If $\{y_j\}_j$ is a $\|\quad \|_2$-dense sequence in $(N_0)_1$ and
$F_n$ is an increasing sequence of finite subsets of $G$ with
$\cup_n F_n = G$, then by the contradiction assumption and the
definition of $\Cal X_0$, there exists $\varepsilon_0 > 0$ such
that for each $n$ there exists $x_n \in (N_0)_1$ with $\|[y_j,
x_n]\|_2 \leq 1/n$, $\forall j \leq n$, $\|\hat{x}_n -
p_{F_n}(\hat{x}_n)\|_2 \geq \varepsilon_0$ and $\|
e_n(\hat{x}_n)-x_n\|_2 \leq 1/n$. (Like in Section 1, if $m \geq
1$ then $e_m$ denotes the orthogonal projection of $L^2(M)$ onto
$L^2(\Sigma_{g \in G} B_m u_g)$.)

But then $x=(x_n)_n \in (N_0^\omega)_1$ satisfies $[x, N_0]=0$, $x
\in \Cal X$ and $\|x - E_{N^\omega \rtimes G}(x)\|_2 \geq
\varepsilon_0$. Since the first two conditions imply $x \in \Cal
X_0$, the third condition contradicts Theorem 3.2. \hfill Q.E.D.

\proclaim{3.4. Corollary} $(\Cal X_0' \cap \langle M^\omega,
N^\omega \rangle)_+$ contains non-zero finite projections.
\endproclaim
\vskip .05in \noindent {\it Proof}. By Lemma 3.3, there exists $F
\subset G$ finite such that $\|p_{F}(\hat{v})\|_2 \geq 1/2$,
$\forall v \in \Cal U(\Cal X_0)$. Note that $p_F=\Sigma_{g\in F}
u_g e_{N^\omega} u_g^* \in \langle M^\omega, e_{N^\omega}\rangle$.
Thus, if $\|\quad \|_{2, Tr}$ denotes the canonical trace on
$\langle M^\omega, e_{N^\omega}\rangle$ (see for instance Sec. 2
in [Po4]) then we have:
$$
Tr (p_F vp_F v^*) = \Sigma_{h, k \in F}
\|E_{N^\omega}(u_h^*vu_k)\|_2^2
$$
$$
= \Sigma_{h,k \in F} \| E_{N^\omega} (vu_{kh^{-1}})\|_2^2 \geq
\|p_F(v)\|_2^2 \geq 1/4.
$$

Thus, if we take $a \in \langle M^\omega, e_{N^\omega} \rangle$ to
be the unique element of minimal norm $\|\quad \|_{2,Tr}$ in
$\overline{\text{\rm co}} \{vp_{F}v^* \mid v \in \Cal U(\Cal
X_0)\}$ then $0 \leq a \leq 1$, $Tr(a) \leq |F| < \infty$ and
$Tr(p_F a) \geq 1/4$. This implies $a\neq 0$ while by uniqueness
we have $vav^*=a$, $\forall v\in \Cal U(\Cal X_0)$. Thus, any
spectral projection corresponding to an interval $[c, 1]$ with $c
> 0$ small enough will do. \hfill Q.E.D.

\vskip .05in

The most important technical result of this section shows that if
in addition to the general assumptions set forth at the beginning
of this section (i.e. $\sigma$ clustering, $pMp = M_0$, $pL(G)p
\supset L(G_0)$) we also assume $pL(G)p \supset L(G_0)$ and $N_0$
abelian, then the von Neumann algebra $\Cal X_0$ is ``large'' in
$M^\omega$. More generally we have:

\proclaim{3.5. Theorem} If $N_0$ is abelian and there exists an
infinite subgroup $G_0' \subset G_0$ such that $pL(G)p \supset
L(G'_0)$ then $\Cal X_0' \cap pM^\omega p = N_0^\omega$.
\endproclaim
\vskip .05in \noindent {\it Proof}. Assume that $x=(x_n)_n \in
\Cal X_0' \cap pM^\omega p$ but $x \not\in N_0^\omega$. Since
$N_0$ is abelian, $N_0'\cap N_0^\omega=N_0^\omega$ so $\Cal
X_0=\Cal X \cap N_0^\omega$ is a von Neumann subalgebra of
$N_0^\omega$ (cf. 3.2) and $N_0^\omega \subset \Cal X_0'\cap
pM^\omega p$.

Thus, by replacing $x$ by $x-E_{N_0^\omega}(x)$, we may assume $x
\perp N_0^\omega$. We may further take $x_n = \Sigma_g y_{n,g}
u^0_g$ with $y_{n,g} \in N_0, y_{n,e}=0$, $\forall n$. Moreover,
by Kaplansky's density theorem we may assume $y_{n,g}=0$ for all
$g \in G_0 \setminus F_n$, for some finite subsets $F_n \subset
G_0$.

At this point we need the following immediate consequences of the
results in Section 2:

\proclaim{3.6. Lemma} If $\sigma$ is clustering and $N_0$ is
abelian, then for any sequence $\{h_n\}_n \subset G'_0$ that tends
to infinity in $G'_0$ and any $y_0 \in N_0$, we have $\tilde{y}_0=
(\sigma_0({h_n})(y_0))_n \in \Cal X_0$.
\endproclaim
\vskip .05in \noindent {\it Proof}. The sequence $(u^0_{h_n})_n$
lies in $L(G)$ and tends weakly to $0$, so $U = (u^0_{h_n})_n \in
L(G)^\omega \ominus L(G)$. By Lemma 2.4 we get $Uy_0 \in \Cal X$
and then by 2.1.2$^\circ$ we get $(Uy_0)U^* \in \Cal X$. Since
$N_0$ is abelian, it follows that $Uy_0U^* \in N_0'\cap N_0^\omega
= N_0^\omega$ as well. Thus $\tilde{y}_0=Uy_0U^* \in \Cal X_0$.
\hfill Q.E.D.

\proclaim{3.7. Lemma} If $\sigma$ is clustering and $N_0$ is
abelian, then ${\sigma_0}_{|G_0'}$ is $n$-mixing, $\forall n\geq
1$.
\endproclaim
\vskip .05in \noindent {\it Proof}. It is sufficient to prove that
if $y_0, y_1, ..., y_n \in N_0$ and $\{h_n^i\}_n\subset G'_0$ are
sequences tending to $\infty$, $1 \leq i \leq n$, then $\underset
n \rightarrow \infty \to \lim \tau(y_0 \Pi_{1\leq i \leq n}
\sigma_0({h^i_n})(y_i)) = 0$. This amounts to showing that
$U_i=(u^0_{h^i_n})_n \in M^\omega$ satisfy $\tau(y_0 Y)= \tau(y_0)
\tau(Y)$, where $Y=\Pi_{1\leq i \leq n}U_iy_iU^*_i$.

As in the proof of 3.6, $h^i_n \rightarrow \infty$ implies $U_i
\in L(G)^\omega \ominus L(G)$. Applying 3.6 and 3.2, we get $Y \in
\Cal X_0$. But $E_M(\Cal X)\subset \Bbb C1$, thus $E_M(Y) =
\tau(Y)1$ and finally
$$
\tau(y_0 Y)= \tau(E_M(y_0Y))=\tau(y_0E_M(Y))=\tau(y_0)\tau(Y). $$
\hfill Q.E.D.

\vskip .05in \noindent {\it Proof of} 3.5 ({\it continuation}).
Since $\sigma_0$ satisfies condition $(1.2.1)$ for $n=2$, by
applying 1.3 and $(1.2.1)$ to sequences of the form
$(\sigma_0({h_n})(y_0))_n$ for $y_0\in N_0$, $h_n \in G'_0, h_n
\rightarrow \infty$ (which by Lemma 3.6 belong to $\Cal X_0$), it
follows that for any finite sets $Y\subset N_0, F \subset G_0$,
any $m, k \geq 1$ and any $\varepsilon>  0$ there exists $n
> k$ such that

$$
|\tau (y_0\sigma_0(h_{n})(y_1)\sigma_0(gh_{n})(y_2)) - \tau(y_0)
\tau(\sigma_0(h_{n})(y_1)\sigma_0(gh_{n})(y_2))| \leq \varepsilon
\tag 3.5.1
$$
for any $y_0, y_1, y_2 \in Y$, $g \in F$.

Let $q \in \Cal P(N_0)$ with $\tau(q)=1/2$ and let $F' \subset
G'_0$ finite such that $\tau(q \sigma_0(g')(q))\leq 9/32$ for $g'
\in G'_0 \setminus F'$ (this is possible because
${\sigma_0}_{G_0'}$ is 2-mixing, hence mixing). Thus,
$\tau(q\sigma_0(g')(q)) + \tau((1-q)\sigma_0(g')(1-q)) \leq 9/16$.
Since $\sigma_0$ is a free action, there exists a finite partition
$\{q_i\}_i \subset \Cal P(N_0)$ refining $\{q, 1-q\}$ such that
$\Sigma_i \tau(q_i \sigma_0(g')(q_i)) \leq 9/16$, for all $g' \in
F', g'\neq e$. Thus, since
$$
\Sigma_i \|q_i \sigma_0(g')(q_i)\|_2^2 = \Sigma_i \tau(q_i
\sigma_0(g')(q_i))
$$
$$
\leq \tau(q\sigma_0(g')(q)) + \tau((1-q)\sigma_0(g')(1-q))
$$
$$
=
\|q\sigma_0(g')(q)\|_2^2 + \|(1-q)\sigma_0(g')(1-q)\|_2^2,
$$
altogether we get

$$
\Sigma_i \|q_i \sigma_0(g')(q_i)\|_2^2 \leq 9/16, \forall g'\in
G_0\setminus \{e\} \tag 3.5.2
$$

Assume $k_1 < k_2 < ... < k_{n-1}$ are given. By applying
$(3.5.1)$ for $\varepsilon = 2^{-4}\|q_i\|_2^2 \|y_{n,g}\|_2^2$,
it follows that there exists $k_n > k_{n-1}$ large enough such
that for all $g \in F_n$ and all $y_{n,g} \neq 0, g\in F_n,$ we
have:
$$
\| (\sigma_0(h_{k_n})(q_i) \sigma_0(gh_{k_n})(q_i)) y_{n,g}\|_2^2
\tag 3.5.3
$$
$$
\leq  (\|\sigma_0(h_{k_n})(q_i) \sigma_0(gh_{k_n})(q_i)\|_2^2 +
2^{-4}\|q_i\|_2^2) \|y_{n,g}\|_2^2
$$

Since the supports $F_n$ of the elements $x_n$ do not contain $e$
and since $g \in F_n$ implies $hgh^{-1} \neq e$ for all $h\in
G_0$, by applying first $(3.5.3)$ and then  $(3.5.2)$ we have the
estimates:
$$
\|x_n\|_2^2 = \|\Sigma_i(\sigma_0(h_{k_n})(q_i) x_n
\sigma_0(h_{k_n})(q_i))\|_2^2
$$
$$
=\Sigma_i \Sigma_{g \in F_n} \| \sigma_0(h_{k_n})(q_i)    y_{n,g}
u^0_g \sigma_0(h_{k_n})(q_i)\|_2^2
$$
$$
=\Sigma_i \Sigma_{g\in F_n} \| (\sigma_0(h_{k_n})(q_i)
\sigma_0(gh_{k_n})(q_i)) y_{n,g}\|_2^2
$$
$$
\leq  \Sigma_i \Sigma_{g \in F_n} (\|\sigma_0(h_{k_n})(q_i)
\sigma_0(gh_{k_n})(q_i)\|_2^2 + 2^{-4}) \|y_{n,g}\|_2^2
$$
$$
=  \Sigma_{g \in F_n} (\Sigma_i \|q_i
\sigma_0(h_{k_n}^{-1}gh_{k_n})(q_i)\|_2^2 + 2^{-4}\|q_i\|_2^2)
\|y_{n,g}\|_2^2
$$
$$
\leq (10/16) \Sigma_{g \in F_n} \|y_{n,g}\|_2^2 =(5/8)
\|x_n\|_2^2,
$$
a contradiction.

\hfill Q.E.D.

\vskip .05in \noindent {\bf 3.8. Remark.} Note that in the above
we have actually proved the following more general result: Let
$\sigma_0$ be a free action of an infinite discrete group $G_0$ on
a finite von Neumann algebra $(N_0, \tau_0)$ and denote $M_0 = N_0
\rtimes_{\sigma_0} G_0$. Assume $\sigma_0$ is {\it weakly
$2$-mixing}, i.e. $\forall F \subset N_0$ finite $\exists h_n
\rightarrow \infty$ in $G_0$ such that
$$
\underset n,m  \rightarrow \infty \to \lim |\tau (y_0
\sigma_0({h_n})(y_1)
\sigma_0({h_m})(y_2))-\tau(y_0)\tau(\sigma_0({h_n})(y_1)
\sigma_0({h_m})(y_2))| =0, \forall y_i\in F.
$$
Let $x\in M_0^\omega$ be so that $E_{N_0^\omega}(x)=0$. Then given
any $\varepsilon > 0$ there exists a partition $\{q_i\}_i$ of $1$
with projections in $N_0$ and a sequence $h_n \rightarrow \infty$
in $G_0$ such that the partition
$\tilde{q}_i=(\sigma_0(h_n)(q_i))_n \in \Cal P(N_0^\omega)$
satisfies $\|\Sigma_i \tilde{q}_i x \tilde{q}_i\|_2 \leq
\varepsilon$.

\heading 4. Cartan conjugacy criteria
\endheading

We now use Section 3 to prove a criterion for unitary (or inner)
conjugacy of Cartan subalgebras. It shows that if $M$ can be
realized in two ways as a group measure space construction,
$M=L^\infty(X, \mu) \rtimes_{\sigma} G = L^\infty(X_0, \mu_0)
\rtimes_{\sigma_0} G_0$, with $\sigma$ clustering, then the
unitary conjugacy of the group von Neumann algebras $L(G), L(G_0)$
entails the unitary conjugacy of the Cartan subalgebras
$L^\infty(X, \mu)$, $L^\infty(X_0, \mu_0)$. In fact, in its most
general form (see 4.2), this statement assumes only a ``corner''
of $M$ to have a second group measure space decomposition $A_0
\rtimes_{\sigma_0} G_0$, with merely $L(G'_0) \subset pL(G)p$ for
some infinite subgroup $G_0' \subset G_0$.

Combined with the result in ([Po4]), showing that if $G$ is ICC,
$G_0$ is w-rigid and $\sigma$ is malleable then the unitary
conjugacy of $L(G_0)$ into $L(G)$ holds true automatically, this
result already allows us to show that modulo perturbation by inner
automorphisms any $\theta : M_0 \simeq M$ must take the Cartan
decompositions $L^\infty(X_0, \mu_0) \subset M_0$, $L^\infty(X,
\mu)\subset M$ onto each other.

We fix the following general assumptions/notations, throughout the
rest of the paper: \vskip .05in \noindent {\bf 4.1. Notations}.
Let $\sigma : G \rightarrow {\text{\rm Aut}} (A, \tau)$, $\sigma_0
: G_0 \rightarrow {\text{\rm Aut}} (A_0, \tau)$ be free, ergodic
actions of countable discrete groups $G, G_0$ on abelian von
Neumann algebra $(A, \tau)$, $(A_0,\tau)$. Denote $M=A
\rtimes_\sigma G$, $M_0 = A_0 \rtimes_{\sigma_0} G_0$ and
$L(G)\subset M$ (respectively $L(G_0)\subset M_0$) the von Neumann
algebra generated by the canonical unitaries $\{u_g\}_{g\in
G}\subset M$ (resp. $\{u_h^0\}_{h\in G_0}\subset M_0$)
implementing the action $\sigma$ (resp. $\sigma_0$).

\vskip .05in

\proclaim{4.2. Theorem} With the notations $4.1$, assume that
$\sigma$ is clustering and there exists a projection $p \in L(G)$
such that $pMp= M_0$, $pL(G)p \supset L(G'_0)$, for some infinite
subgroup $G_0'\subset G_0$. Then there exists a partial isometry
$v \in M$ such that $v^*v=p$, $vv^* \in A$ and $vA_0v^* = Avv^*$.
\endproclaim
\vskip .05in \noindent {\it Proof}. Let $\Cal X_0 \subset
pM^\omega p$ be defined as in Section 3. By its definition, $\Cal
X_0$ is contained in $A_0^\omega$ and by Theorem 3.2, $\Cal X_0$
is a von Neumann algebra. Also, by Theorem 3.5 we have $\Cal
X_0'\cap pM^\omega p = A_0^\omega$.

Due to Corollary 3.4, we are in the position of applying the
``intertwining theorem'' (A.1 in [Po2]). It follows that there
exist non-zero projections $q=(q_n)_n \in A^\omega$,
$p^0=(p^0_n)_n \in A_0^\omega \subset pM^\omega p$ and a partial
isometry $w=(w_n)_n \in M^\omega$ such that $q_n \in A, p_n^0 \in
A_0$, $w_n^* w_n = p^0_n, w_nw_n^* = q_n$,
$\tau(q_n)=\tau(p_n^0)=\tau(q)$, $\forall n$, and $wA_0^\omega
w^*= A^\omega q$. This implies that for $\varepsilon = \tau(q)/2$
and large enough $n$ we have $w_nA_0w_n^* \subset_\varepsilon
Ap_n^0$. Indeed, because if for each $n$ we could find $u_n \in
\Cal U(w_nA_0w_n^*)$ such that $\|E_{A^\omega} (u_n) - u_n\|_2
\geq \tau (q)/2$, then $u=(u_n)_n$ would satisfy $wuw^* \not\in
A^\omega$, contradicting $wA_0^\omega w^* \subset A^\omega$.

By taking $\tilde{A}_0 \subset M$ to be a Cartan subalgebra of
$M$ such that $p\in \tilde{A}_0$ and $A_0=\tilde{A}_0p$, and
applying (Corollary 1 in [Po5]) to the Cartan subalgebras $A,
\tilde{A}_0$ of $M$, it follows that there exists a partial
isometry $v \in M$ with $v^*v=p$ and $vA_0v^*=Avv^*$. \hfill
Q.E.D.

\vskip .05in \noindent {\bf 4.3. Remarks}. If instead of
(Corollary 1 in [Po5]) we use results from ([PoSiSm]), or some of
Christensen's pioneering results on perturbations of subalgebras
in II$_1$ factors ([Ch]), then the argument at the end of the
proof of 4.2 shows the following: If $A, A_0$ are arbitrary
maximal abelian $^*$-subalgebras in a II$_1$ factor $M$ such that
$A_0^\omega, A^\omega$ are unitary conjugate in $M^\omega$ then
$A_0, A$ are unitary conjugate in $M$.

\vskip .05in

In the next Section, we will prove that if in addition to the
conditions in 4.2 we also assume $G$ ICC, then the situation is in
fact much more rigid: $p$ must equal 1 and $v$, which follows a
unitary, can be taken to satisfy $v\{u^0_h\}_hv^*=\{u_g\}_g$,
modulo the scalars ! In particular, this will show that the
inclusion $L(G_0) \subset pL(G)p$ must be an equality. We prove
here a weaker result in this direction, showing that $A_0, pL(G)p$
follow mutually orthogonal (i.e., $E_{L(G)}(a_0)=0$, $\forall
a_0\in A_0, \tau(a_0)=0$; see [Po7]). This is needed in the proof
of the criteria for conjugacy of actions in the next section (5.1,
5.2).

\proclaim{4.4. Lemma} With the notations $4.1$, assume that for
some projection $p\in L(G)$ we have  $pMp\supset M_0$, $pL(G)p
\supset L(G_0)$ and there exists a partial isometry $v\in M$ such
that $v^*v=p$, $vA_0v^* =Avv^*$. \vskip .05in $(i)$. If $\sigma_0$
is weakly mixing then $A_0$ follows orthogonal to $L(G)$.

\vskip .05in $(ii)$. If $\sigma_0$ is weakly mixing and $\sigma$
is mixing, then $\sigma_0$ follows mixing.
\endproclaim
\vskip .05in \noindent {\it Proof}. $(i)$. Note first that
$\sigma_0$ weakly mixing implies that $\forall y\in A_0$ with
$\tau(y)=0$ there exist $h_n \in G_0$ such that $\sigma_0(h_n)(y)$
tends to $0$ in the weak operator topology, as $n \rightarrow
\infty$.

Assume there exists $y_0 \in A_0 \subset pMp$ such that
$\tau(y_0)=0$, $x_0=E_{L(G)}(y_0) \neq 0$. Let $h_n \in G_0$ be so
that $u^0_{h_n}y_0u_{h_n}^{0*}=\sigma_0(h_n)(y_0)$ tends weakly to
$0$. Expecting on $L(G)$ and taking into account that
$$
E_{L(G)}(u^0_{h_n}y_0u_{h_n}^{0*})=u^0_{h_n}E_{L(G)}(y_0)u_{h_n}^{0*}
=u^0_{h_n}x_0u_{h_n}^{0*},
$$
it follows that the sequence $\tilde{x}_0=(u_{h_{n}}^0 x_0
u^{0*}_{h_{n}})_n \subset L(G)$ is convergent to $0$ in the weak
operator topology. Thus, as an element in the ultrapower algebra
$L(G)^\omega$, $\tilde{x}_0$ belongs to $L(G)^\omega\ominus L(G)$,
while $\|\tilde{x}_0\|_2=\|x_0\|_2 \neq 0$.

But $A_0 = v^*Av$ implies $A_0^\omega = v^*A^\omega v$, and thus
$A_0^\omega$ is contained in $A^\omega \rtimes G$, the von Neumann
algebra generated by $A^\omega$ and the unitaries $u_g \in M$,
$g\in G$, regarded as constant sequences in $M^\omega$. Thus, if
we denote $\tilde{y}_0=(u^0_n y_0 u^{0*}_n)_n \in A_0^\omega$ then
$\tilde{y}_0 \in A^\omega \rtimes G$. By commuting squares it
follows that  $E_{L(G)^\omega}(A^\omega \rtimes_\sigma
G)=L(G)^\omega \cap A^\omega \rtimes_\sigma G \subset L(G)$. Thus
$E_{L(G)^\omega}(\tilde{y}_0) \in L(G)$. On the other hand
$$
E_{L(G)^\omega}(\tilde{y}_0)=(E_{L(G)}(u^0_n y_0 u^{0*}_n))_n =
\tilde{x}_0 \in L(G)^\omega\ominus L(G).
$$
This contradiction ends the proof of part $(i)$.

$(ii)$. Let $h_n \rightarrow \infty$ in $G_0$ and $a_i \in A_0$
with $\tau(a_i)=0$, $i=1,2$. Since by $(i)$ we have $A_0 \perp
L(G)$, it follows that $\tau(a_2)=0$ implies $a_2 \perp L(G)$. But
then by 2.6 we get $\underset n \rightarrow \infty \to \lim
\tau(a_2 u^0_{h_n}a_1 u^{0*}_{h_n}) = 0$, showing that $\sigma$ is
mixing. \hfill Q.E.D.

Our last result in this section shows that under the same
assumptions as in 4.5, the condition $\sigma_0$ mixing is in fact
not much of a restriction:

\proclaim{4.5. Lemma} With the notations $4.1$, assume that for
some projection $p\in L(G)$ we have $pMp \supset M_0$, $pL(G)p
\supset L(G_0)$. Assume also that either $pMp = M_0$ or that there
exists $v\in M$ partial isometry such that $v^*v=p$, $vA_0v^* =
Avv^*$. If $\sigma$ is mixing then there exist a normal subgroup
$G_1 \subset G_0$ of finite index $n$ and a projection $p_1 \in
\Cal P(A_0)$ of trace $\tau(p_1)=\tau(p)/n$ that commutes with
$\{u_h^0\}_{h\in G_1},$ such that if $A_1\overset \text{\rm def}
\to =A_0p_1$, $\sigma_1(h)\overset \text{\rm def} \to =
\sigma_0(h)_{|A_1}$, $u^1_h\overset \text{\rm def} \to=u^0_hp_1,
h\in G_1,$ then: \vskip .05in $(i)$. $\sigma_1$ is free and mixing
on $A_1$ and $M_1\overset \text{\rm def} \to = p_1M_0p_1=A_1 \vee
{\{u^1_h\}_{h\in G_1}}''$ is isomorphic to $A_1 \rtimes_{\sigma_1}
G_1$, with $\{u^1_h\}_h$ as canonical unitaries implementing
$\sigma_1$. \vskip .05in $(ii)$. $p_1Mp_1 \supset M_1$, with
equality iff $pMp=M_0$; $p_1L(G)p_1 \supset L(G_1) =
{\{u^1_h\}_{h\in G_1}}''$, with equality iff $pL(G)p=L(G_0)$.
Moreover, if $v\in M$ satisfies $v^*v=p$, $vA_0v^*=Avv^*$ then
$v_1=vp_1\in M$ is a partial isometry satisfying $v_1^*v_1=p_1$,
$v_1A_1v_1^* = Av_1v_1^*$.
\endproclaim
\vskip .05in \noindent {\it Proof}. If $B_0\subset A_0$ is a
finite dimensional $\sigma_0$-invariant subalgebra of $A_0$ then
${\text{\rm Sp}}L(G_0)B_0={\text{\rm Sp}}B_0L(G_0)$. Thus, (3.1 in
[Po4]) applied to $Q_0 = L(G_0)p \oplus (1-p)L(G)(1-p) \subset A
\rtimes_\sigma G$ implies that $B_0$ is contained in $L(G)$ and
therefore in $pL(G)p$. This shows that the von Neumann subalgebra
$B\subset A_0$ generated by all such $B_0 \subset A_0$ is
contained in $pL(G)p$ as well. If $B$ would be infinite
dimensional then by the ergodicity of $\sigma$ it follows diffuse.
Since the normalizer of $B$ in $pMp$ contains both $A_0$ and
$\{u^0_h \mid h\in G_0\}$, which together generate $M_0$, by
applying again (3.1 in [Po4]) it follows that $M_0 \subset
pL(G)p$. This gives an obvious contradiction when we assume
$pMp=M_0$. If we assume instead that $A_0$ is conjugate to a
corner of $A$ by some partial isometry $v\in M$, then $A_0$ would
contain an orthonormal family of unitary elements that are
perpendicular to $A$ (the powers of any generating Haar unitary
for $B$ would do), contradicting (A.1 in [Po2], or 2.1 in [Po4]).

Thus $B$ must be finite dimensional and by the ergodicity of
$\sigma_0$ all its minimal projections have the same trace. Let
$n=\text{\rm dim}(B)$ and choose $p_1\in \Cal P(B)$ a minimal
projection. Let also $G_1=\{h \in G_0 \mid \sigma_0(h)(b)=b,
\forall b\in B\}$. It is immediate to check that all conditions
$(i), (ii)$  are satisfied. \hfill Q.E.D.

\heading 5.  Conjugacy criteria for actions
\endheading

We now prove a very general conjugacy (or isomorphism) criterion
for actions of groups. Thus, we show that under very general
conditions if a group measure space factor $M$ has two crossed
product decompositions $M=A\rtimes_\sigma G = A_0
\rtimes_{\sigma_0} G_0$ in a way that the group algebras $L(G),
L(G_0)$ coincide, then modulo conjugacy by an inner automorphism
of $M$ the two decompositions are identical ! The result is in
fact surprisingly more general, requiring only the inclusion of
one group algebra into the other and allowing amplifications.

We recall that if $K \vartriangleleft G$ is a normal subgroup and
$\kappa \in \text{\rm Char}_G(K)$ is an Ad$(G)$ invariant
character of $K$ then one associates to it the scalar 2-cocycle
$\mu_\kappa \in \text{\rm H}^2(G/K,\Bbb T)$ as follows: Denote by
$g \mapsto g'$ a lifting of $G/K$ into $G$ then for $g_1, g_2 \in
G/K$ put $\mu_\kappa(g_1,g_2)\overset \text{\rm def} \to =
\kappa(g_1' g_2' ((g_1g_2)')^{-1})$. It is easy to verify that
$\mu_\kappa$ doesn't depend on the lifting and that it defines
indeed a 2-cocycle on $G/K$.

\proclaim{5.1. Theorem} With the general notations $4.1$, assume
$\sigma$ is clustering and there exists a projection $p \in L(G)$
such that $pMp= M_0$, $pL(G)p \supset L(G_0)$. Then $p$ is central
in $L(G)$, $\tau(p)^{-1}$ is an integer, $pL(G)p = L(G_0)$ and
there exist a normal subgroup $K \subset G$ with
$|K|=\tau(p)^{-1}$, an ${\text{\rm Ad}}(G)$-invariant character
$\kappa$ of $K$, with trivial $\mu_{\kappa} \in H^2(G/K)$, and a
unitary element $u\in pL(G)p$ such that:

\vskip .05in $(i)$. $p=|K|^{-1} \Sigma_{k \in K} \kappa(k) u_k$
and $u_kp=\kappa(k) p, \forall k\in K$.

\vskip .05in $(ii)$. $uA_0u^* = A^Kp$, where $A^K=\{a\in A \mid
\sigma_k(a)=a, \forall k\in K\}$.

\vskip .05in $(iii)$. $\{u_gp \mid g\in G\} = \{uu^0_hu^* \mid
h\in G_0\}$ modulo multiplication by scalars. More precisely,
there exist an isomorphism $\delta: G_0 \simeq G/K$, a lifting
$G/K \ni g \mapsto g' \in G$ and a map $\alpha : G_0 \simeq G/K
\rightarrow \Bbb T$ such that $\text{\rm Ad}(u)(u^0_h)=
\overline{\alpha(h)} u_{\delta(h)'} p$, $\forall h\in G_0$, and
$\partial\alpha =\mu_\kappa$.

If in addition $G$ has no non-trivial finite normal subgroups
(e.g. $G$ ICC or torsion free) then $K=\{e\}$, $p=1$ and $\delta$
is an isomorphism $\delta:G_0\simeq G$.
\endproclaim
\vskip .05in \noindent {\it Proof}. Since ICC groups have no
normal subgroups other than $K=\{e\}$, the last part follows
trivially from $(i)$ and $(iii)$.

By Theorem 4.2, the conditions in the hypothesis imply the
existence of a partial isometry $v \in M$ with $v^*v=p$ and
$vA_0v^* = Avv^*$. We are thus reduced to proving $5.1.(i)-(iii)$
when this last condition is added to the list of assumptions.
Notice also that by 3.7 the clustering condition on $\sigma$
implies $\sigma_0$ mixing. We in fact prove the following more
general result:

\proclaim{5.2. Theorem} With the general notations $4.1$, assume
there exists a projection $p \in L(G)$ such that:

\vskip .05in $(a)$. $p=1_{M_0}$, $pMp\supset M_0$, $pL(G)p \supset
L(G_0)$.

\vskip .05in  $(b)$. One of the following holds true: $p=1$,
$\sigma_0$ weakly mixing; or $p$ arbitrary, $\sigma_0$ weakly
mixing, $\sigma$ mixing.

\vskip .05in  $(c)$. There exists $v\in M$ such that $v^*v = p$
and $vA_0{v}^* = Avv^*$.

\vskip .05in

Then $\tau(p)^{-1}$ is an integer and there exist a subgroup $G'
\subset G$, a normal subgroup $K \subset G'$ with
$|K|=\tau(p)^{-1}$, a character $\gamma$ of $G_0$, an ${\text{\rm
Ad}}(G')$-invariant character $\kappa$ of $K$, with trivial
$\mu_{\kappa} \in H^2(G'/K)$, and a unitary element $u\in pL(G)p$
such that $p'=upu^*$ is central in $L(G')=\{u_g \mid g\in G'\}''$
and we have:

\vskip .05in  $(i)$. $p'=|K|^{-1} \Sigma_{k \in K} \kappa(k) u_k
\in \Cal Z(L(G'))$, $u_kp'=\kappa(k) p', \forall k\in K$.

\vskip .05in  $(ii)$. $uA_0u^* = A^Kp'=p'Ap'$, where $A^K=\{a\in A
\mid \sigma_k(a)=a, \forall k\in K\}$.

\vskip .05in  $(iii)$. $\{u_gp' \mid g\in G'\} = \{uu^0_hu^* \mid
h\in G_0\}$ modulo multiplication by scalars. More precisely,
there exist an isomorphism $\delta: G_0 \simeq G'/K$, a lifting
$G'/K \ni g \mapsto g' \in G'$ and $\alpha: G_0 \simeq G'/K
\rightarrow \Bbb T$ such that $\text{\rm Ad}(u)(u^0_h)=
\overline{\alpha(h)} u_{\delta(h)'} p'$, $\forall h\in G_0$, and
$\partial \alpha = \mu_\kappa$.

\vskip .05in If in addition $pMp = M_0$ and $G$ has no non-trivial
finite normal subgroups (e.g. $G$ is ICC or torsion free), then
condition $(b)$ is redundant, $K=\{e\}$, $p=1$, $G'=G$ and
$\delta$ is an isomorphism $\delta:G_0\simeq G$.
\endproclaim
\vskip .05in \noindent {\it Proof}. Let us first prove the last
part of the statement, assuming we have proved $(i)-(iii)$ under
conditions $(a)-(c)$. Thus, we are under the additional
assumptions $pMp=M_0$ and $G$ is ICC. By Lemma 4.5, there exists a
subgroup $G_1 \subset G_0$ of index $n < \infty$ and a projection
$p_1\in \Cal P(A_0)$ of trace $\tau(p_1)=\tau(p)/n$ that commutes
with $\{u^0_h\}_{h\in G_1}$ such that if we denote $A_1=A_0p_1$,
$\sigma_1(h)=\sigma_0(h)_{|A_1}$ then $\sigma_1$ is mixing. But
then conditions $(a)-(c)$ are satisfied for $\sigma_1, G_1, p_1$,
$v_1=vp_1$ instead of $\sigma_0, G_0, p, v$. Thus, by $(i)-(iii)$
there exist $G'_1\subset G$, with a finite normal subgroup $K_1
\subset G_1'$ satisfying $|K_1| = \tau(p_1)^{-1}$, and a unitary
element $u_1 \in p_1L(G)p_1$ such that $p_1'=u_1p_1u_1^*$ is
central in $L(G_1')$ and $u_1A_1u_1^* = A^{K_1}p_1'=p_1'Ap_1'$.
Thus, after conjugation by $u_1$ we may assume $p_1$ satisfies
$A_0=A^{K_1}p_1=p_1Ap_1$ and $L(G_1)= p_1L(G'_1)p_1=L(G'_1)p_1$.
But sp$L(G_1) A_1$ is dense in $M_1$ and $p_1(A \rtimes
G'_1)p_1=\overline{\text{\rm sp}} A^{K_1} L(G_1')p_1$. Since
$p_1Mp_1=M_1$, by expecting these equalities on $L(G)$ we get
$L(G_1')p_1=p_1L(G)p_1$, which in particular implies $[G:G'] <
\infty$. But then $G$ ICC implies $G_1'$ ICC (indeed, for if
$G_1'$ would have a non-trivial element $h$ with finite conjugacy
class in $G_1'$ then $h$ would also have finite conjugacy class in
$G$, due to $[G:G_1'] < \infty$). Thus $K_1=\{e\}$ forcing
$p_1=1$.  From $L(G_1')p_1=p_1L(G)p_1$ we also get $G_1'=G$ and
$\delta$ follows an isomorphism because $K_1=\{e\}$.

The proof of $(i)-(iii)$ will take the rest of this Section, as
well as Section 6. The following observation will be used several
times:

\vskip .05in
\noindent
$(5.2.0)$. If $u'$ is a partial isometry in
$L(G)$ with ${u'}^*u'=p$ and we substitute $p$ by $u'{u'}^*$,
$L(G_0)$ by $u'L(G_0){u'}^*$, $M_0$ by $u'M_0{u'}^*$, $v$ by
$v{u'}^*$ and $A_0$ by $u'A_0{u'}^*$ then conditions $(a)-(c)$ are
still satisfied. \vskip .05in

We first prove the result assuming the Fourier expansion
$v=\Sigma_g a_gu_g$ is finitely supported, with all $a_g \in A$
having finite spectrum. Under this assumption the proof
simplifies, allowing the ideas to become more transparent. The
proof of the general case is postponed to the next section.

Since all $a_g$ have finite spectrum and there are only finitely
many of them, there exists $q \neq 0$ in $\Cal P(A)$, $q \leq
vv^*$, such that $qa_g = c_g q$, $\forall g$, for some scalars
$c_g \in \Bbb C$. Thus, if we denote $w=\Sigma_g c_gu_g$, then
$w\in L(G)$ and $qw=qv$. Since $v^*(Aq)v\subset A_0$ it follows
that $w^*qAqw \subset A_0$. But $A$ is perpendicular to $L(G)$ and
by Lemma 4.4 $A_0$ is also perpendicular to $pL(G)p$.
Consequently,

$$
\tau(q)/\tau(p) p = E_{L(G)}(v^*qv)
=E_{L(G)}(w^*qw)=w^*E_{L(G)}(q)w =\tau(q) w^*w.
$$

Hence, $u'= \tau(p)^{-1/2} w$ is a partial isometry in $L(G)$ with
${u'}^*u'=p$. Moreover, it satisfies
$qv{u'}^*=qw{u'}^*=\tau(p)^{-1/2} u'{u'}^*$. Thus, by the remark
at the beginning of the proof it follows that in addition to
conditions $5.2.(a)-(c)$ we may also assume:

$$
qv=\tau(p)^{-1/2} qp \tag 5.2.1
$$

We will prove that in fact conditions $5.2.(a)-(c)$ together with
$(5.2.1)$ imply $5.2.(i)-(iii)$ for $u=p$.

Since $u_h^0$ normalizes $A_0$ and $vA_0v^* = Avv^*$, $vu_h^0v^*$
is in the normalizing groupoid of $A$. Thus $qvu_h^0v^*q$ is also
in the normalizing groupoid of $A$, implying that
$qvu_h^0v^*q=\Sigma_g a_g^hu_g$ for some partial isometries $a_g^h
\in A$. On the other hand $u_h^0 \in L(G_0) \subset pL(G)p$, so
$u_h^0=\Sigma_g c^h_g u_g$ for some scalars $c^h_g$. Together with
$qv=\tau(p)^{-1/2}qp$ this implies:
$$
\Sigma_g a^h_g u_g=qvu_h^0v^*q=\tau(p)^{-1} qu_h^0 q =
\tau(p)^{-1} \Sigma_g c^h_g q\sigma_g(q) u_g. \tag 5.2.2
$$

Identifying  the coefficients of the most left and right Fourier
expansions in $(5.2.2)$, it follows that $a_g^h=\tau(p)^{-1} c_g^h
q\sigma_g(q)$, $\forall g$. In particular, if we denote
$F^h=\{g\in G \mid  a^h_g \neq 0\}=\{g \in G \mid c^h_g \neq 0,
\sigma_g(q)q \neq 0\}$, then $g \in F^h$ implies
$1=\|a_g^h\|=\tau(p)^{-1} |c_g^h|$, and thus $|c_g^h|=\tau(p)$.
Since $\tau(p)=\|u_h^0\|_2^2=\Sigma_g |c^h_g|^2$, this also shows
that $|F^h| \leq \tau(p)^{-1}$. Summarizing:
$$
|c^h_g|=\tau(p), \forall g\in F^h; \Sigma_g |c^h_g|^2 = \tau(p);
|F^h| \leq \tau(p)^{-1}. \tag 5.2.3
$$

When $p=1$ then $qv=q$, $Aq=A_0q$, and $(5.2.3)$ shows that if
$F^h \neq \emptyset$ then $F^h$ is a single point set $\{g(h)\}$,
with $c^h=c^h_{g(h)}$ satisfying $|c^h|=1$ and $u^0_h = c^h
u_{g(h)}$. Denote by $G_1\subset G_0$ the subgroup of elements $h
\in G_0$ with $u_h^0$ a scalar multiple of some $u_{g(h)}$. If
$h_0 \in G_0$ is given, then by the weak mixing property of
$\sigma_0$ there exists $h_n \in G_0$ such that $\underset n
\rightarrow \infty \to \lim \tau(\sigma_0(h_n)(q)q)=\tau(q)^2$ and
$\underset n \rightarrow \infty \to \lim
\tau(\sigma_0(h_nh_0)(q)q)=\tau(q)^2$. Noticing that for $h'\in
G_0$ we have $F^{h'}\neq \emptyset$ $\Leftrightarrow$ $qu_{h'}^0q
\neq 0$ $\Leftrightarrow$ $\sigma_0(h')(q)q\neq 0$, it follows
that $h_n, h_nh_0 \in G_1$, which in turn implies $h \in G_1$.
This shows that $G_0=G_1$. Finally notice that $Aq=A_0q$ implies
$$
A_0\sigma_0(h)(q)=u^0_h(A_0q){u^0_h}^* = u_{g(h)}(Aq)u_{g(h)}^*=
A\sigma_{g(h)}(q) \tag 5.2.4
$$
and that by the ergodicity of $\sigma_0$ we have $1=\vee_h
\sigma_0(h)(q)=\vee_h \sigma_{g(h)}(q)$. This shows that $A=A_0$
and letting $G'=\{g(h) \mid h\in G_0\}$, $5.2.(i)-(iii)$ follow
trivially.

Assume now $\sigma$ is mixing but $p$ is arbitrary. Notice that if
for some $h$ we have $|F^h|=\tau(p)^{-1}$ then $(5.2.3)$ forces
$F^h$ to be equal to the support set $S^h=\{g \in G \mid c^h_g
\neq 0\}$ of the Fourier expansion $u^0_h = \Sigma_g c^h_g u_g$,
and $|c^h_g|=\tau(p), \forall g\in S^h$. We will prove that this
is in fact the case for $\forall h \in G_0$, but first show it is
true for all $h$ outside a finite subset of $G_0$.

By the mixing property, there exists $E\subset G$ finite such that
$\tau(\sigma_g(q)q) \geq \tau(q)^2/2$, $\forall g \in G\setminus
E$. Define $\varepsilon_0 > 0$ to be the distance from
$\tau(p)^{-1}$ to the closest integer $\neq \tau(p)^{-1}$. Since
$u^0_h$ tends to $0$ in the weak operator topology, as $h
\rightarrow \infty$, there exists $E_0 \subset G_0$ finite such
that $|c^h_g|^2=|\tau(u^0_hu_g^*)|^2 < \varepsilon_0
\tau(p)^2/(2|E|+2)$, $\forall g \in E$, $\forall h \in
G_0\setminus E_0$. Thus, if $h \in G_0 \setminus E_0$ and $g \in
G\setminus E$ then $g \in S^h$ $\Leftrightarrow$ $g \in F^h$
$\Leftrightarrow$ $|c^h_g|=\tau(p)$ (the last equivalence by
$(5.2.3)$).

If for some $h \in G_0 \setminus E_0$ we have $|F^h| <
\tau(p)^{-1}$, then from the above we get
$$
\tau(p)=\Sigma_g |c_g^h|^2 = \Sigma_{g \in E} |c^h_g|^2 +
\Sigma_{g \in F^h} |c_g^h|^2 \tag 5.2.5
$$
$$
=\Sigma_{g \in E} |c^h_g|^2 + |F^h| \tau(p)^2 < \varepsilon_0
\tau(p)^2/2 + |F^h| \tau(p)^2
$$
$$
=(\varepsilon_0/2 + |F^h|) \tau(p)^2 \leq (\tau(p)^{-1} -
\varepsilon_0/2) \tau(p)^2 = \tau(p)-\varepsilon_0\tau(p)^2/2,
$$
a contradiction. We have thus shown that $\tau(p)^{-1}$ is an
integer and that for $h \in G_0 \setminus E_0$ we have: $S^h$ has
$\tau(p)^{-1}$ elements; $|c^g_h| = \tau(p), \forall g\in S^h$.
Conclusions $(i)$ and $(iii)$ in 5.2 follow then from the
following general fact:

\proclaim{5.3. Lemma} Let $G$ be a discrete group, $p\in
L(G)={\{u_g\}_g}''$ a projection with $\tau(p)=1/n$ for some
integer $n \geq 1$. Let $\Cal W \subset pL(G)p$ denote the set of
unitary elements $w \in pL(G)p$ such that $S^w=F^w$, where $S^w =
\{g \in G \mid c^w_g \neq 0\}$, $F^w=\{ g \in G \mid
|c^w_g|=\tau(p) \}$, $c^w_g$ being the coefficients of the Fourier
expansion $w=\Sigma_g c^w_g u_g$. Assume $\Cal W \neq \emptyset$.
Then $\Cal W$ is a subgroup of $\Cal U(pL(G)p)$ with neutral
element equal to $p$ and there exist a subgroup $G' \subset G$, a
normal subgroup $K \subset G'$ with $|K|=n$ and an ${\text{\rm
Ad}}(G')$-invariant character $\kappa$ of $K$ such that $p=n^{-1}
\Sigma_{k\in K} \kappa(k) u_k$, $u_k p= \kappa(k)u_k$, $\forall
k\in K$, and $\Cal W = \{\alpha u_t p \mid t\in G', \alpha \in
\Bbb T \}$. Moreover, if $G'/K \ni h \mapsto h' \in G'$ is a
lifting, then $\{u_{h'}p \mid h \in G'/K\}$ is a copy of the
projective left regular representation of $G'/K$ with scalar
$2$-cocycle $\mu_\kappa \in \text{\rm H}^2(G'/K)$.
\endproclaim
\vskip .05in \noindent {\it Proof}. Since $\Cal W \cdot \Cal W
\subset \Cal U(pL(G)p)$ and any $u\in\Cal U(pL(G)p)$ satisfies
$\tau(u^*u)=\tau(p)$, it follows that for all $w,w'\in \Cal W$ we
have:
$$
|F^w|=n; \Sigma_g |c^w_g|^2=\tau(p)=\Sigma_t |c^{ww'}_t|^2;
\Sigma_g |c^w_g|=1. \tag 5.3.1
$$

We also have
$$
ww'=\Sigma_{t \in S^{ww'}} (\Sigma_{g \in S^w}
c_g^w c^{w'}_{g^{-1}t})u_t = \Sigma_t c^{ww'}_t u_t,
$$
and there
are at most $n^2$ non-zero elements $c_t^{ww'}, t\in G,$ and they
satisfy the  trivial inequalities
$$
|c^{ww'}_t| \leq \Sigma_g |c^w_g| |c^{w'}_{g^{-1}t}| \leq \tau(p)
\tag 5.3.2
$$
$$
\Sigma_t |c^{ww'}_t| \leq \Sigma_{g,g'} |c^w_g| |c^{w'}_{g'}|=1
\tag 5.3.3
$$

Moreover, the first inequality in $(5.3.2)$ becomes equality iff
$c^w_g c^{w'}_{g^{-1}t}, g\in G,$ are all equal, when non-zero,
while the second inequality in $(5.3.2)$ becomes equality iff $g
\mapsto g^{-1}t$ is a bijection from $S^w$ to $S^{w'}$.

Notice that the maximum of the expression $\Sigma_{i=1}^{n^2}
x_i^2$ over the set of $n^2$-tuples $\{(x_i)_i \mid 0 \leq x_i
\leq 1/n, \forall i, \Sigma_j x_j \leq 1\}$ is equal to $1/n$ and
it is attained when $n$ of the $x_i$'s are equal to $1/n$ and all
others are equal to $0$. Indeed, this is because
$$
x_i^2 + x_j^2 \leq (\text{min}\{x_i+x_j, 1/n\})^2 +
(\text{max}\{0, x_i+x_j-1/n\})^2,
$$
with equality if and only if $\{x_i,x_j\}=\{\text{min}\{x_i+x_j,
1/n \}, \text{max}\{0, x_i+x_j-1/n\}\}$.

Applying this to $(x_i)_i = (|c_g^{ww'}|)_g$, it follows that
$\Sigma_t |c^{ww'}_t|^2 \leq 1/n = \tau(p)$ with equality iff
there are exactly $n$ non-zero elements $c^{ww'}_t$ and they all
have absolute value $\tau(p)$, i.e. iff $S^{ww'}=F^{ww'}$. Since
$(5.3.1)$ shows that the equality does hold true, it follows that
$ww'\in \Cal W$. In addition, the conditions for equality in
$(5.3.2)$ are verified, i.e.
$$
S^{w'}=(S^w)^*t, c^w_g c^{w'}_{g^{-1}t}=c^{ww'}_t \tau(p), \forall
w,w'\in \Cal W, t\in S^{ww'}, g\in S^w \tag 5.3.4
$$

Since we clearly have $\Cal W^* = \Cal W$, this shows that
$p=ww^*$ belongs to $\Cal W$ and that $\Cal W$ is a group with $p$
as neutral element. Denote $K = S^p$ and note that $K=K^*$. By
$(5.3.4)$ applied to $w=p=w'$, it follows that $K$ is a group with
$n$ elements. Since $p$ is a projection of the form $p=|K|^{-1}
\Sigma_{k \in K} \alpha_k u_k$, with $\alpha_k=c^p_k/\tau(p)$
having absolute value 1, $\forall k\in K$, $\kappa(k) \overset
\text{\rm def} \to = \alpha_k, k\in K$, gives a character of $K$.

Also, by $(5.3.4)$ applied first to $w=p$, then to $w'=p$, it
follows that for all $w\in \Cal W$ there exist $t, t' \in G$ with
$S^w=Kt=t'K$. This implies $KtK=tK$, so $t$ normalizes $K$. Thus,
if we denote $G'=\{t\in G \mid \exists w\in \Cal W, S^w=Kt\}$ then
$G'$ follows a group with $K\subset G'$ a normal subgroup.
Moreover, $w$ is a scalar multiple of $pu_t$ and $[u_t,p]=0,
\forall t\in G'$. In particular, the character $\kappa$ of $K$
must be ${\text{\rm Ad}}(G')$-invariant and it satisfies
$u_kp=\kappa(k)p, \forall k\in K$. The last part of the statement
is now trivial, since $\tau(u_g p)=0$, $\forall g \in G' \setminus
K$ and the 2-cocycle $\mu_\kappa \in \text{\rm H}^2(G'/K)$
satisfies $u_{(h_1)'} p' u_{(h_2)'} p'=\mu_\kappa(h_1,h_2)
u_{(h_1h_2)'} p'$ by the way it is defined. \hfill Q.E.D. \vskip
.05in

We can now end the proof of $5.2$ under the extra assumption that
the Fourier expansion $v=\Sigma_g a_g u_g$ is finitely supported,
with all $a_g$ having finite spectrum. Thus, by the formula
$p=\Sigma_{k\in K} \kappa(k) u_k$ it follows that $A^K=A \cap
\{p\}'$ with $p$ a Jones projection implementing the
$\tau$-preserving conditional expectation $E_{A^K}$ of $A$ onto
$A^K$. Since $qv=\tau(p)^{-1/2}qp$, the projection $q\in A$
satisfies $pqp=E_{A^K}(q)p = \tau(p)q'p$ for some unique
projection $q'\in A^K$. By $(5.2.1)$ we thus get
$$
A_0q'p=A_0(pqp)=A_0(v^*qv) = v^*(Aq)v = p(Aq)p = A^Kq'p \tag 5.4.1
$$
The equality $A_0p=A^Kp$ in $5.2.(ii)$ follows now from this
identity, the ergodicity of $\sigma_0$ on $A_0$ and the inclusion
$\{u^0_h\}_{h\in G_0} \subset \{\alpha u_gp \mid g\in G', \alpha
\in \Bbb T\}$, the same way we deduced the similar equality in the
case $p=1$ in $(5.2.4)$. Part $5.2.(iii)$ is an immediate
consequence of the last part of Lemma 5.3.

\heading 6. End of proof of 5.2
\endheading

We now prove 5.2 in general, without any extra-assumption on the
intertwiner $v=\Sigma_g a_g u_g$. The idea is to prove first an
approximate version of condition $(5.2.1)$, and then re-do steps
$(5.2.2)-(5.2.5)$ using it in-lieu of $(5.2.1)$.

\proclaim{6.1. Lemma} Let $(A,\tau), G, \sigma, M=A \rtimes_\sigma
G, \{u_g\}_g$ be as in $4.1$. Let $p \in \Cal P(L(G))$ and assume
$v\in M$ is a partial isometry with $v^*v=p$, $vv^* \in A$ and
$v^*Av$ perpendicular to $L(G)$. Then there exists a partial
isometry $u' \in L(G)$ and a decreasing sequence of non-zero
projections $q_n \in Avv^*$ such that ${u'}^*u'=p$ and
$$
\|q v  - \tau(p)^{-1/2} qu'\|^2_2 < 2^{-n} \tau(q), \forall q\in
\Cal P(Aq_n) \tag 6.1
$$
\endproclaim
\vskip .05in \noindent {\it Proof}. Let $v = \Sigma_g a_g u_g$ be
the Fourier expansion of $v$ in $M = A \rtimes_\sigma G$ and
denote $f_0=vv^* \in \Cal P(A)$. Since the infinite sum $\Sigma_g
a_g a_g^*$ is increasing and tends to $f_0$ in the strong operator
topology, or equivalently in $\|\cdot \|_1$, there exists $f \neq
0$ in $\Cal P(Af_0)$ such that $\Sigma_g a_g a_g^*f$ is convergent
to $f$ in the uniform norm. Thus, there exists an increasing
sequence of finite sets $F_n \subset G$ such that $\cup_n F_n = G$
and
$$
\|\Sigma_{g \not\in F_n} a_ga_g^*f\| \leq 2^{-n-3} \tag 6.1.1
$$

By the Gelfand-Neimark theorem $Af=C(\Omega)$ for some compact set
$\Omega$. Choose $\omega_0 \in \Omega$ and denote $c_g =
a_g(\omega_0), g\in G$. Thus $\Sigma_g c_g\overline{c_g} = 1$,
implying that $w = \Sigma_g c_g u_g$ belongs to
$L^2(L(G))=\ell^2(G)$. By the continuity of the functions $a_g \in
C(\Omega)$ at $\omega_0$, we can choose recursively $q_n \in \Cal
P(Af)$ (corresponding to open-closed neighborhoods of $\omega_0$
in $\Omega$) such that $q_n \leq q_{n-1}$ and
$$
\|a_g q_n - c_g q_n \| \leq 2^{-n-2}|F_n|^{-1/2}, \forall g \in
F_n \tag 6.1.2
$$
Using first the Cauchy-Schwartz inequality and then $(6.1.2)$ we
get for any $q\in \Cal P(Aq_n)$:
$$
|\Sigma_{g\in F_n} \tau(\overline{c_g} (a_g-c_g)q)| =
|\Sigma_{g\in F_n} \tau(\overline{c_g}q (a_g-c_g)q)| \tag 6.1.3
$$
$$
\leq (\Sigma_{g\in F_n} \|c_g q \|_2^2)^{1/2} (\Sigma_{g\in F_n}
\|(a_g-c_g) q \|_2^2)^{1/2}
$$
$$
\leq (\Sigma_{g \in F_n} |c_g|^2 \tau(q))^{1/2}(\Sigma_{g \in F_n}
2^{-2n-4}|F_n|^{-1} \tau(q))^{1/2}
$$
$$
\leq (\Sigma_{g \in G} |c_g|^2 \tau(q))^{1/2} 2^{-n-2}
\tau(q)^{1/2} = 2^{-n-2} \tau (q).
$$

Similarly, by first using the Cauchy-Schwartz inequality and then
$(6.1.1)$ we have:

$$
|\tau(\Sigma_{g\not\in F_n} \overline{c_g}a_g q)| \tag 6.1.4
$$
$$
\leq (\Sigma_{g\not\in F_n} |c_g|^2 \tau(q))^{1/2}
(\tau(\Sigma_{g\not\in F_n} a_ga_g^*q))^{1/2} \leq 2^{-n-3}
\tau(q).
$$

Using $(6.1.3)$ and $(6.1.4)$ we finally obtain the estimates:
$$
\|qw-qv\|_2^2 = 2 \tau(q) - 2 {\text{\rm Re}} \tau(qvw^*q) \tag
6.1.5
$$
$$
=2 \tau(q) - 2 {\text{\rm Re}} \tau(\Sigma_g \overline{c_g}a_g q)
= 2 \tau(q) -2 {\text{\rm Re}} \tau(\Sigma_{g\in F_n}
\overline{c_g}a_g q) - 2 {\text{\rm Re}} \tau(\Sigma_{g\not\in
F_n} \overline{c_g}a_g q)
$$
$$
= 2 \tau(q) - 2\Sigma_{g \in F_n} c_g \overline{c_g} \tau(q) +
2\Sigma_{g\in F_n} \tau(\overline{c_g} (a_g-c_g)q) - 2 {\text{\rm
Re}} \tau(\Sigma_{g\not\in F_n} \overline{c_g}a_g q)
$$
$$
\leq 2 \tau (q) - 2\tau(q) + 2 \cdot 2^{-n-3} \tau(q) + 2\cdot
2^{-n-2} \tau(q) + 2 \cdot 2^{-n-3} \tau(q) = 2^{-n}\tau(q).
$$
This shows that if we put $u'=\tau(p)^{1/2}w$ and $q_n$ are chosen
as above, then the inequality $(6.1)$ is satisfied $\forall q\in
\Cal P(Aq_n)$.

We still have to prove that ${u'}^*u'=p$. To this end, note that
$w\in L^2(L(G))$ and $q_n \in Aq_n \perp L(G)$ implies
$E_{L(G)}(w^*q_nw)=w^*E_{L(G)}(q_n)w=\tau(q_n)w^*w$ (these
calculations make sense in $L^1(L(G))$). On the other hand,
$v^*q_nv \in v^*(Aq_n)v \perp L(G)$ by hypothesis, so that we also
have $E_{L(G)}(v^*q_nv)=(\tau(q_n)/\tau(p)) p$. But by the
Cauchy-Schwartz inequality and $(6.1)$ we have
$$
\|E_{L(G)}(w^*q_nw-v^*q_nv)\|_1  \leq \|w^*q_nw-v^*q_nv\|_1
$$
$$
\leq \|w^*q_n w - w^* q_n v\|_1 + \|w^* q_n v - v^* q_n v\|_1
$$
$$
\leq \|w^*q_n\|_2 \|q_n w - q_n v\|_2 + \|w^* q_n  - v^* q_n\|_2
\|q_n v\|_2
$$
$$
= 2 \|q_n w - q_n v\|_2 \|q_n\|_2 \leq  2^{-n+1} \tau(q_n).
$$
Thus, for all $n \geq 1$ we have
$$
\|\tau(q_n)/\tau(p)) p - \tau(q_n) w^*w \|_1 =
\|E_{L(G)}(w^*q_nw-v^*q_nv)\|_1 \leq 2^{-n+1} \tau(q_n)
$$
implying that $u'=\tau(p)^{1/2} w$ satisfies $\|{u'}^*u'-p\|_1
\leq \tau(p) 2^{-n+1}$, $\forall n$. Thus ${u'}^*u'=p$.

\hfill Q.E.D.

\vskip .1in Let now $u'\in L(G)$ and $q_n \in \Cal P(Avv^*)$ be as
given by Lemma 6.1. Since for $q\in \Cal P(Aq_n)$ we have
$\|qv{u'}^* - \tau(p)^{-1/2}qp\|_2 = \|qv-\tau(p)^{-1/2}qu'\|_2
\leq 2^{-n}\|q\|_2$, by the remark at the beginning of the proof
of 5.2 in Section 5, we may assume there exists a decreasing
sequence of non-zero projections $q_n \in Avv^*$ such that in
addition to the assumptions $5.2.(a)-(c)$ we also have:
$$
\|qv-\tau(p)^{-1/2} qp\|_2 \leq 2^{-n}\|q\|_2, \forall q\in \Cal
P(Aq_n). \tag 6.2
$$

We will prove that in fact conditions $5.2.(a)-(c)$ together with
$(6.2)$ imply $5.2.(i)-(iii)$ for $u=p$. In particular, when $p=1$
this amounts to showing that $\{u^0_h\}_h \subset \{ \alpha u_g
\mid g\in G, \alpha \in \Bbb T\}$ and $A=A_0$.

As in Section 5, let $u^0_h = \Sigma_g c^h_g u_g$ be the Fourier
expansion of $u^0_h \in L(G_0) \subset pL(G)p$, with $c^h_g \in
\Bbb C$. Since $u^0_h$ are unitaries in $pL(G)p$, we have
$\tau(p)=\tau(u_h^0u^{0*}_h)=\Sigma_g |c^h_g|^2$, $\forall h\in
G_0$. Note that $q_n u^0_h q_n = \Sigma_g c^h_g \sigma_g(q_n)q_n
u_g$. We also consider the Fourier expansion $q_n v u^0_h v^* q_n
= \Sigma_g a^{h,n}_g u_g$. Since the left term is in the
normalizing groupoid of $Aq_n \subset q_n M q_n$, it follows that
$a^{h,n}_g, g\in G,$ (resp. $\sigma_{g^{-1}}(a^{h,n}_g)$, $g\in
G$) are partial isometries in $Aq_n$ with mutually orthogonal
supports. But unfortunately we no longer have the equality $q_nv
u^0_h v^* q_n=\tau(p)^{-1} q_nu^0_h q_n$ that would make possible
the identification of the Fourier coefficients of the two sides,
as in $(5.2.2)$. However, an ``approximation'' of that equality
does hold true due to $(6.2)$, and this will be good enough for us
to prove the following version of $(5.2.3)$:

\proclaim{6.3. Lemma}  Let $F^h_n = \{g\in G \mid  a^{h,n}_g \neq
0\}$. For all $g \in F^h_n$ we have: \vskip .05in $(i)$. $\|a^h_g
-c^h_g \tau(p)^{-1} \sigma_g(q_n)q_n \| \leq 2^{-n+1}$ and if $n
\geq 2$ then $|a^h_g| =\sigma_g(q_n)q_n$. \vskip .05in $(ii)$.
$|1-\tau(p)^{-1} |c^h_{g}| | \leq 2^{-n+1}$ and $|F^h_n| \leq
\tau(p)^{-1} + 2^{-n+2}$. \vskip .05in $(iii)$. $|c^h_g| \leq
\tau(p)$ and $F^h_n \supset F^h_{n+1}$, $\forall n$.
\endproclaim
\vskip .05in \noindent {\it Proof}. If $F^h_n = \emptyset$ then
there is nothing to prove. If  $F^h_n$ is non-empty then fix $g_0
\in F^h_n$. Denote $\overline{q}=q_n \sigma_{g_0}(q_n)$ and note
that $|a^{h,n}_{g_0}| \leq \overline{q}$. For each projection $q
\in \Cal P(A\overline{q})$ denote $q'=\sigma_{g_0^{-1}}(q) \leq
q_n$. Note that if $q'\leq \overline{q}$ then $q u_{g_0} q' =
qu_{g_0}$ and if $q\leq |a^{h,n}_{g_0}|$ then $qvu^0_h v^*
q'=qa^{h,n}_{g_0}u_{g_0}$.

By $(6.1)$ applied to $q,q' \leq q_n$ we have $\|qv -
\tau(p)^{-1/2} qp\|_2 \leq 2^{-n}\|q\|_2$ and $\|v^*q' -
\tau(p)^{-1/2} pq'\|_2 \leq 2^{-n}\|q'\|_2$. Using this and the
inequality $\|xy\|_2 \leq \|x\| \|y\|_2$ we get:

$$
\|qvu^0_h v^* q' - \tau(p)^{-1} qu^0_h q' \|_2 \tag 6.3.1
$$
$$
\leq \|qvu^0_h v^* q' - \tau(p)^{-1/2}qu^0_h v^*q'\|_2 +
\|\tau(p)^{-1/2}qu^0_h v^*q'-\tau(p)^{-1} qu^0_h q' \|_2
$$
$$
\leq \|qv  - \tau(p)^{-1/2}qp \|_2 \|u^0_hv^*q'\| + \tau(p)^{-1/2}
\| qu^0_h\| \| v^*q'-\tau(p)^{-1/2} pq' \|_2
$$
$$
\leq 2^{-n} \|q\|_2 + 2^{-n} \|q'\|_2 = 2^{-n+1} \|q\|_2.
$$

For any $d$ in the spectrum of $a^{h,n}_{g_0}$ in $A\overline{q}$
and any $\varepsilon
> 0$, take $0\neq q\in \Cal P(A\overline{q})$ to be a spectral projection
of $a^{h,n}_{g_0}$ in $A\overline{q}$ such that $\|a^{h,n}_{g_0} q
- d\| \leq \varepsilon$. Using $(6.3.1)$ we then get:
$$
2^{-n+1} \|q\|_2 \geq \|qvu^0_h v^* q' - \tau(p)^{-1} qu^0_h q'
\|_2
$$
$$
= \|(a^{h,n}_{g_0}q - \tau(p)^{-1} c^h_{g_0} q) u_{g_0} -
\Sigma_{g \neq g_0} c^h_g q \sigma_g(q') u_g \|_2
$$
$$
\geq \|a^{h,n}_{g_0}q - \tau(p)^{-1} c^h_{g_0}q\|_2 \geq \|d q -
\tau(p)^{-1} c^h_{g_0} q\|_2 - \varepsilon \|q\|_2
$$
$$
= |d - \tau(p)^{-1} c^h_{g_0}| \|q\|_2 - \varepsilon \|q\|_2.
$$
Since $d\in {\text{\rm spec}}(a^{h,n}_{g_0})$ and $\varepsilon >
0$ were arbitrary, $\| a^{h,n}_{g_0} - \tau(p)^{-1}
c^h_{g_0}\overline{q} \|\leq 2^{-n+1}$. This also shows that if $n
\geq 2$ then $|a^{h,n}_{g_0}|=\overline{q}$. Moreover,
$|1-\tau(p)^{-1} |c^{h}_{g_0}| | \leq \|a^{h,n}_{g_0} -
\tau(p)^{-1} c^h_{g_0} \overline{q}\| \leq 2^{-n+1}$.

From the equality $\Sigma_{g\in G} |c^h_g|^2 = \tau(p)$ and the
inequalities $|c^h_g| \geq \tau(p)(1-2^{-n+1})$ we get
$$
\tau(p) \geq \Sigma_{g \in F^h_n} |c^h_g|^2 \geq |F^h_n| \tau(p)^2
(1-2^{-n+1})^2,
$$
showing that $|F^h_n| \leq \tau(p)^{-1} (1-2^{-n+1})^{-2} \leq
\tau(p)^{-1} (1+2^{-n+2})$.

The sequence of sets $F^h_n$ follows decreasing in $n$ by the
definition of $F^h_n$ and by the fact that the sequence $q_n$ is
decreasing. Also, by the Cauchy-Schwartz inequality we have
$|c^h_g| = |\tau (u^0_gu_g^*)| \leq \|u^0_g\|_2 \|u_gp\|_2
=\tau(p)$. \hfill Q.E.D.

\vskip .1in \noindent {\it Proof of} 5.2 {\it in case} $p=1,$
$\sigma_0$ {\it weakly mixing.} Since $p=1$, condition $(6.2)$
becomes:
$$
\|qv-q\|_2 \leq 2^{-n} \|q\|_2, \forall q\in \Cal P(Aq_n) \tag
6.4.1
$$
for some decreasing sequence of non-zero projections $q_n$ in
$Avv^*$. By the weak mixing property of $\sigma_0$, there exists a
sequence $\{h_m\}_m \in G_0$ such that
$$
\underset m \rightarrow \infty \to \lim
\tau(\sigma_0(h_mh)(v^*q_nv)v^*q_nv)=\tau(q_n)^2, \forall n \geq
1, h\in G_0. \tag 6.4.2
$$

Note that for $h'\in G_0$ and $n \geq 1$ we have $F^{h'}_n \neq
\emptyset \Leftrightarrow q_n vu^0_{h'}v^* q_n \neq 0
\Leftrightarrow \sigma_0(h')(v^*q_nv) v^*q_nv \neq 0$. Let now $h
\in G_0$ be an arbitrary element and $n \geq 3$. Taking $h'=h_m$
and then $h'=h_mh$ in these equivalences, with $m$ large enough,
it follows by $(6.4.2)$ that $F^{h_m}_{n}, F^{h_mh}_{n} \neq
\emptyset$. By Lemma 6.3 we have $|F^{h_mh}_n|, |F^{h_m}_n| \leq
1+2^{-n+2}\leq 3/2$. Thus $F^{h_m}_{n}=\{g\}$,
$F^{h_mh}_{n}=\{g'\}$, for some $g=g(n,h_m), g'=g(n,h_mh)$.

Moreover, by 6.3 again, there exist unit scalars $c,c' \in \Bbb T$
such that $|c - c^{h_n}_g| \leq 2^{-n+2}$, $|c' - c^{h_nh}_{g'}|
\leq 2^{-n+2}$. Recalling that $u^0_{h_n}=\Sigma_g c^{h_n}_g u_g$,
we get
$$
\|u^0_{h_n} - cu_g\|^2_2 = |c^{h_n}_g - c|^2 + (1 - |c^{h_n}_g|^2)
$$
$$
\leq 2^{-2n+4} + 2^{-n+2}(1+ 2^{-n+2}) \leq 2^{-n+3}.
$$
Similarly $\|u^0_{h_n} - cu_g\|^2_2 \leq 2^{-n+3}$. Since $u_h^0=
u_{h_n^{-1}}^0 u^0_{h_n h}$, by the triangle inequality we deduce:
$$
\|u_h^0 - c^{-1}c' u_{g^{-1}g'}\|_2
$$
$$
\leq \|u_{h_n^{-1}}^0 u^0_{h_n h}-u_{h_n^{-1}}^0 (c'u_{g'})\|_2 +
\|u_{h_n^{-1}}^0 (c'u_{g'}) - (c^{-1}u_{g^{-1}})(c'u_{g'}\|_2
$$
$$
=\|u^0_{h_n h}-c'u_{g'}\|_2 +\|u_{h_n^{-1}}^0 -
c^{-1}u_{g^{-1}}\|_2 \leq 2^{(-n+5)/2}.
$$

Since $n \geq 3$ can be taken arbitrarily large in these
inequalities, independently of $h \in G_0$, this shows that
$\{u^0_h\}_h \subset \{\alpha u_g \mid g \in G, \alpha \in \Bbb T
\}$. Since $\{u^0_h\}_h$ is a representation of the group $G_0$,
there exist a character $\gamma$ of $G_0$ and an isomorphism
$\delta:G_0 \simeq G'\subset G$ such that $u^0_{h}=\gamma(h)
u_{\delta(h)}$, $\forall h\in G_0$.

This proves all but $(ii)$ in 5.2. Proving $(ii)$ amounts to show
that under the extra-assumption $(6.2)$, $A$ and $A_0$ follow
equal. But by $(6.2)$, under the ``infinitesimal'' projections
$q_n\in A$, $v^*q_nv\in A_0$ the two algebras ``almost'' coincide.
Since $\{u^0_h\}_h$ is contained in $\{u_g\}_g$ modulo the scalars
and they  act ergodically on $A_0$ resp. $A$, this local ``almost
coincidence'' extends to a global ``almost coincidence'' of $A$,
$A_0$, which due to the arbitrariness of the approximations shows
that $A=A_0$. The full details of this argument are carried out in
$(6.6.1)-(6.6.3)$, where the case ``$p$ arbitrary'' is covered.

\hfill Q.E.D.

\vskip .1in \noindent {\it Proof of} 5.2 {\it in the case}
$\sigma_0$ {\it weakly mixing,} $\sigma$ {\it mixing.} By
$4.4.(i)$, since $\sigma$ is mixing and $\sigma_0$ weakly mixing,
it follows that $\sigma_0$ is mixing as well. Let $\varepsilon_0 >
0$ be the distance from $\tau(p)^{-1}$ to the closest integer
$\neq \tau(p)^{-1}$ and take $0 < \varepsilon \leq \varepsilon_0
(1+\tau(p)^{-1})^{-1}/2$.

Since $\sigma$ is mixing there exists a finite subset $E \subset
G$ such that
$$
|\tau(\sigma_g(q_n)q_n) -\tau(q_n)^2| \leq \varepsilon
\tau(q_n)^2, \forall g \in G \setminus E. \tag 6.5.1
$$
Similarly, since $\sigma_0$ is mixing and
$\sigma_0(h)(y)=u^0_hyu^{0*}_h$, $\forall y \in A_0=Ap, h \in
G_0$, by taking into account that $v^*q_nv \in A_0$ and
$\tau_{M_0}(y) =\tau(p)^{-1} \tau(y)$, $\forall y \in A_0=Ap$, it
follows that there exists $E_0 \subset G_0$ finite such that
$$
|\tau(u^{0}_hv^*q_nvu^{0*}_h v^*q_nv)-\tau(p)^{-1}\tau(q_n)^2|
\leq \varepsilon \tau(q_n)^2, \forall h \in G_0\setminus E_0. \tag
6.5.2
$$
Moreover, since $\{u^0_h\}_h$ tends weakly to zero as $h
\rightarrow \infty$, it follows that we may also assume $E_0$ is
large enough to ensure
$$
|c^h_g|=|\tau(u^0_hu_g^*)| \leq \varepsilon/|E|, \forall g \in E,
\forall h\in G_0 \setminus E_0 \tag 6.5.3
$$

But $q_nvu^0_hv^*q_n = \Sigma_g a^{h,n}_gu_g$, so using the fact
that $\|\Sigma_g a^{h,n}_gu_g\|_2^2=\Sigma_g \tau(|a^{h,n}_g|)$
and $\|q_nvu^0_hv^*q_n\|_2^2=\tau(u^{0}_hv^*q_nvu^{0*}_h
v^*q_nv)$, by $(6.5.2)$ it follows that
$$
|\Sigma_g \tau(|a^{h,n}_g|)-\tau(p)^{-1}\tau(q_n)^2| \leq
\varepsilon \tau(q_n)^2, \forall h \in G_0\setminus E_0. \tag
6.5.4
$$
On the other hand, by Lemma 6.3, for each $a^{h,n}_g \neq 0$ we
have $\tau(|a^{h,n}_g|)=\tau(\sigma_g(q_n)q_n)$. Thus, since there
are $|F^{h}_n|$-many such non-zero elements, by first using
$(6.5.4)$ and then $(6.5.1)$ we get
$$
| |F^h_n|-\tau(p)^{-1}| | \tau(q_n)^2=
|F^h_n|\tau(q_n)^2-\tau(p)^{-1}\tau(q_n)^2|
$$
$$
\leq | |F^h_n|\tau(q_n)^2 - \Sigma_g \tau(\sigma_g(q_n)q_n)| +
\varepsilon \tau(q_n)^2
$$
$$
\leq (|F^h_n|+1)\varepsilon \tau(q_n)^2 \leq
(\tau(p)^{-1}+1)\varepsilon \tau(q_n)^2 \leq (\varepsilon_0/2)
\tau(q_n)^2.
$$
Thus, $\tau(p)^{-1}$ is $\varepsilon_0/2$ close to the integer
$|F^h_n|$. By the choice of $\varepsilon_0$, this forces
$\tau(p)^{-1}=|F^h_n|$.

We have thus proved that $\tau(p)^{-1}$ is an integer and that for
each $n$ there exists $E_0 \subset G_0$ finite such that for all
$h \in G_0 \setminus E_0$ we have $|F^h_n|=\tau(p)^{-1}$, while by
$6.3.(ii)$ we also have $|\tau(p)-|c^h_g|| \leq 2^{-n+1}\tau(p)$,
$\forall g \in F^h_n$. Define $\tilde{G}=\{(u_{g_n})_n \mid g_n
\in G\} \subset \Cal U(M^\omega)$ and notice that $L(\tilde{G})$
is naturally isomorphic to the von Neumann subalgebra of
$M^\omega$ generated by $\tilde{G}$. Apply Lemma 5.3 to the set
$\Cal W$ of unitary elements in $pL(\tilde{G})p$ which have
Fourier expansion with $\tau(p)^{-1}$ non-zero coefficients, all
of same absolute value $\tau(p)$. From the above remarks it
follows that if $\{h_n\}_n \subset G_0$ satisfies $h_n \rightarrow
\infty$ then the class of $(u^0_{h_n})_n$ in $M^\omega$ lies in
$\Cal W$. Since for all $h\in G_0$ we have $hh_n \rightarrow
\infty$ as well, it follows that $(u^0_{hh_n})_n \in \Cal W$,
$\forall h\in G_0$. But by 5.3 $\Cal W$ is a group, so the
constant sequence $u^0_h=(u^0_{hh_n})_n (u^0_{h_n})_n^*$ lies in
$\Cal W$, $\forall h\in G_0$.

We have thus shown that the Fourier expansion of $u^0_h$ has
exactly $\tau(p)^{-1}$ elements, all of absolute value $=\tau(p)$,
$\forall h\in G_0$. Thus, if we take $u=p$ then $(i), (iii)$ of
$5.2$ are consequences of Lemma 5.3. In other words, there exist a
subgroup $G'\subset G$, with a normal subgroup $K \subset G'$ with
$\tau(p)^{-1}$ elements and an Ad$(G')$-invariant character
$\kappa$ of $K$ such that $p=\Sigma_k \kappa(k)u_k$ and
$\{u^0_h\}_h=\{u_gp \mid g \in G'\}$, modulo scalars.

In particular, this implies that $p$ is a Jones projection for the
inclusion $A^K \subset A$, i.e. $pap=E_{A^K}(a)p, \forall a\in A$,
and there exists a projection $p_0 \in A$ of trace
$\tau(p_0)=\tau(p)$ satisfying $pp_0p=\tau(p) p$, $A^Kp_0=Ap_0$.

To finish the proof of 5.2 we still need to show that
$A_0=A^Kp=pAp$. To this end note that $(6.2)$ implies the
following ``approximation'' of the equality $A_0(v^*qv)= A^Kq'p$
in $(5.4.1)$:
$$
\|v^*qv - \tau(p)^{-1}pqp\|_2 \leq 2^{-n+1} \|q\|_2, \forall q\in
\Cal P(Aq_n). \tag 6.6.1
$$
Moreover, since in particular $pq_np$ is $2^{-n+1}$-close to a
$\tau(p)$-multiple of a projection, for $n$ large enough it
follows that $pqp$ must be equal to $\tau(p)q'$ where
$q'=\tau(p)^{-1} E_{A^K}(q)$ is a projection in $A^K$, $\forall
q\in \Cal P(Aq_n)$. Thus, with the above notations we may assume
$q_n \leq p_0$ and
$$
\|v^*qv - q'p\|_2 \leq 2^{-n+1} \|q\|_2, \forall q\in \Cal
P(Aq_n). \tag 6.6.1'
$$

Let $\Cal F=\{(e_i, h_i)\}_{i\in I} \subset \Cal P(Aq_n) \times
G_0$ be a maximal family (with respect to inclusion) with the
property that $\{\sigma_{\delta(h_i)}(e_i')\}_i$, resp.
$\{\sigma_0(h_i)(v^*e_iv)\}_i$, are mutually orthogonal. If $f_0 =
\Sigma_i \sigma_0(h_i)(v^*e_iv) \neq p$ (equivalently $f =\Sigma_i
\sigma_{\delta(h_i)}(e_i') \neq 1$), then by the ergodicity of
$\sigma_0$ there exists $0\neq e_0 \in \Cal P(Aq_n)\subset A^Kp_0$
and $h \in G_0$ such that $\sigma_0(h)(v^*e_0v) \leq p-f_0$. By
the maximality of the family $\Cal F$, it follows that
$\sigma_{\delta(h)}(e_0) f \neq 0$. By shrinking $e_0$ if
necessary, we may thus assume $\sigma_{\delta(h)}(e_0) \leq
\sigma_{\delta(h_i)}(e_i')$ for some $i\in I$. Thus, there exists
$e^0_i \leq e_i$ in $\Cal P(Aq_n)$ such that
$$
\sigma_{\delta(h_i)}(e^0_i)=\sigma_{\delta(h)}(e_0),
\sigma_0(h)(v^*e_0v) \sigma_0(h_i)(v^*e^0_iv)=0. \tag 6.6.2
$$

But since for $a \in A^K, a_0 \in A_0$ we have
$\sigma_{\delta(h)}(a)p=u_{\delta(h)} (ap) u^*_{\delta(h)}=u^0_h a
u^{0*}_h$ and $\sigma_0(h)(a_0)=u^0_h a_0 u^{0*}_h=u_{\delta(h)}
a_0 u^*_{\delta(h)}$, by $(6.2)$ we get
$$
\|\sigma_0(h')(v^*qv)-\sigma_{\delta(h')}(q)\|_2 \leq 2^{-n+1},
\forall q \in \Cal P(Aq_n), \forall h'\in G_0. \tag 6.6.3
$$
From $(6.6.2)$ and $(6.6.3)$ we finally get:

$$
2^{1/2} \|e_0\|_2
=\|\sigma_0(h)(v^*e_0v)-\sigma_0(h_i)(v^*e^0_iv)\|_2
$$
$$
\leq \|\sigma_0(h)(v^*e_0v)-\sigma_{\delta(h)}(e_0)\|_2 +
\|\sigma_{\delta(h_i)}(e^0_i) -\sigma_0(h_i)(v^*e^0_iv)\|_2
$$
$$
\leq 2^{-n+2} \|e_0\|_2
$$
which for $n \geq 2$ gives a contradiction. The maximal family
$\Cal F$ must therefore satisfy $\Sigma_i
\sigma_0(h_i)(v^*e_iv)=p$,  $\Sigma_i \sigma_{\delta(h_i)}(e_i')
=1.$ Since by $(6.2)$  we have $\|v^*e_i - e_i\|_2 \leq 2^{-n}
\|e_i\|_2$, if we define $V = \Sigma_i u_{h_i}^0 e_iv
u_{h_i}^{0*}$ then clearly $V^*V=VV^*=p$, $V^*(A^Kp)V=A_0$ and
$$
\|V-p\|_2^2 = \Sigma_i \|u_{h_i}^0 e_iv u_{h_i}^{0*} - u_{h_i}^0
e_iu_{h_i}^{0*}\|_2^2 \leq \Sigma_i 2^{-n} \|e_i\|_2^2 = 2^{-n}.
$$
But since $n$ was arbitrary, this shows that $A^Kp=A_0$.

\hfill Q.E.D.

\heading 7. Strong rigidity results
\endheading

The results we prove in this section will typically assume a
``weak equivalence'' between two free ergodic m.p. actions
$\sigma: G \rightarrow {\text{\rm Aut}}(X_0, \mu_0)$, $\sigma: G
\rightarrow {\text{\rm Aut}}(X,\mu)$ satisfying certain conditions
and will derive from that a much ``stronger equivalence''. There
are two types of ``weak equivalence'' assumptions that we will
consider:

{\it Von Neumann equivalence} (abbreviated {\it vNE}) of two
actions $\sigma_0, \sigma$ is an isomorphism of the group measure
space factors $M_0=L^\infty(X_0,\mu_0) \rtimes_{\sigma_0} G_0$,
$M=L^\infty(X,\mu) \rtimes_\sigma G$ associated with $\sigma_0,
\sigma$. Most often we'll in fact consider its ``stabilized
version'', i.e. isomorphism of some ``corners'' of these factors.

{\it Orbit equivalence} (abbreviated {\it OE}) of $\sigma_0,
\sigma$, is an isomorphism $\Delta_0$ of the probability spaces
$(X_0, \mu_0)$, $(X,\mu)$ that takes the equivalence relation
$\Cal R_{\sigma_0}$ given by the orbits of $\sigma_0$ onto the
equivalence relation $\Cal R_\sigma$ given by the orbits of
$\sigma$, almost everywhere. More generally we consider its
``stabilized version'', i.e. an isomorphism $\Delta_0:(X_0,
\mu_0)\simeq (Y,\mu_Y)$, for some $Y\subset X$ with $\mu(Y)\neq
0$, that takes $\Cal R_{\sigma_0}$ onto the equivalence relation
$\Cal R^Y_\sigma$ on $(Y,\mu_Y)$ given by the intersection of the
orbits of $\sigma$ with the set $Y$.

Note that by ([Dy], [FM]) an isomorphism $\Delta_0: (X_0,\mu_0)
\simeq (Y,\mu_Y)$ extends to an algebra isomorphism $\theta_0:M_0
\simeq pMp$, where $p=\chi_Y$, if and only if it takes onto each
other the equivalence relations $\Cal R_{\sigma_0}, \Cal
R^Y_\sigma$, i.e. if it is a (stable) OE. Furthermore, an
isomorphism $\theta:M_0 \simeq pMp$ comes from a (stable) OE iff
it takes the subalgebra of coefficients $L^\infty(X_0,\mu_0)$,
$L^\infty(X,\mu)p$ (also called Cartan subalgebras) onto each
other. Thus, OE $\Rightarrow$ vNE and any OE of actions can be
viewed as a special type of vNE of the actions. Altogether, OE
ergodic theory, which is the study of actions of groups up to
orbit equivalence, is the same as the study of Cartan subalgebra
inclusions $L^\infty(X,\mu) \subset L^\infty(X,\mu)\rtimes_\sigma
G$ underlying the group measure space construction.

{\it Conjugacy of actions} is the typical ``stronger equivalence''
that we'll derive for $\sigma_0, \sigma$. This means an
isomorphism of probability spaces $\Delta : (X_0,\mu_0) \simeq
(X,\mu)$ satisfying $\{\Delta^{-1} \sigma_g \Delta\}_{g\in G}=
\{\sigma_0(h)\}_{h\in G_0}$. In other words, $\sigma(\delta(h))
\Delta = \Delta \sigma_0(h)$, $\forall h\in G_0$, for some
isomorphism of groups $\delta:G_0 \simeq G$. Note that if $\Delta,
\delta$ give a conjugacy of $\sigma_0, \sigma$, then they
implement an isomorphism of the corresponding group measure space
II$_1$ factors $\theta^{\delta,\Delta} :L^\infty(X_0,\mu_0)
\rtimes_{\sigma_0} G_0 \simeq L^\infty(X,\mu) \rtimes_{\sigma} G$
by $\theta^{\delta,\Delta}(\Sigma_h a^0_h u^0_h) = \Sigma_h
\Delta(a_h^0) u_{\delta(h)}$.

In one of the results though (Theorem 7.6), from vNE we'll merely
derive OE, with only certain ``parts'' of $(\sigma_0,G_0),
(\sigma,G)$ conjugate onto each other.

The conditions $\sigma_0, \sigma$ must satisfy are as follows: The
``source'' action $(\sigma_0, G_0)$ has to contain a sufficiently
large subgroup satisfying the relative property (T) (of
Kazhdan-Margulis [Ma]). The ``target'' action $(\sigma,G)$ needs
to have good deformation properties (malleability), examples of
which are the Bernoulli $G$-actions. The precise definitions and
notations, also used in ([Po4]), are as follows:

\vskip .05in
\noindent
{\it 7.0.1. Definition.} A group $G_0$ is
{\it w-rigid}, if it contains an infinite normal subgroup
$H_0\subset G_0$ with the relative property (T) of
Kazhdan-Margulis ([Ma]; see also [dHV]).

\vskip .1in \noindent {\it 7.0.2. Definition} ([Po6]). An infinite
subgroup $\Lambda$ of a group $\Gamma$ is {\it weakly quasi
normal} ({\it wq-normal}) in $\Gamma$ if for any intermediate
subgroup $\Lambda \subset H \varsubsetneq \Gamma$ there exists $g
\in \Gamma \setminus H$ such that $gHg^{-1} \cap H$ infinite. Note
that $\Lambda$ is wq-normal in $\Gamma$ if and only if there
exists an ordinal $\imath$ and a well ordered family of
intermediate subgroups $\Lambda = H_0 \subset H_1 \subset ...
\subset H_{\jmath} \subset ... \subset H_{\imath} = \Gamma$ such
that for each $\jmath < \imath$, $H_{\jmath+1}$ is the group
generated by the elements $g \in \Gamma$ with $gH_{\jmath}g^{-1}
\cap H_{\jmath}$ infinite and such that if $\jmath \leq \imath$
has no ``predecessor'' then $H_{\jmath} = \cup_{n < \jmath} H_n$
(see [Po6]). It is easy to see that there exists a largest group
$\Lambda \subset H'\subset \Gamma$ such that $\Lambda$ is
wq-normal in $H'$, which we will call  the {\it wq-normalizer} of
$\Lambda$ in $\Gamma$.

If $\Lambda=H_0 \subset H_1 \subset ... \subset H_n =\Gamma$ are
all normal inclusions, and $\Lambda$ is infinite, then $\Lambda
\subset \Gamma$ is wq-normal. An inclusion of the form $\Lambda
\subset (\Lambda * K_0) \times K_1$, with $\Lambda, K_1$ infinite
$K_0$ arbitrary, is also wq-normal but cannot be realized as a
sequence of consecutive normal inclusions.

\vskip .1in \noindent {\it 7.0.3. Notation.}  We denote by $w\Cal
T_0$  the class of groups $G_0$ that contain an infinite,
non-virtually abelian subgroup $H$ satisfying: $H\subset G_0$ has
the relative property (T); $H$ is wq-normal in $G_0$.

\vskip .1in \noindent {\it 7.0.4. Notation.} We denote by $w\Cal
T_1$ the class of groups $G_0$ that have a w-rigid, wq-normal
subgroup $H \subset G_0$ such that $\{hgh^{-1} \mid h\in H\}$ is
infinite $\forall g \in G_0, g\neq e$. Note that if $G_0$ is
itself w-rigid ICC then $G_0 \in w\Cal T_1$.
 \vskip .05in
Note that all groups in $w\Cal T_1$ are ICC (by the definition),
while groups in $w\Cal T_0$ are not required to be ICC. The class
$w\Cal T_0$ is clearly closed under normal and finite index
extensions. It is also closed to inductive limits. All infinite
property (T) groups ([Ka]) are both w-rigid and in the class
$w\Cal T_0$, and so are all groups of the form $\Gamma \times K$
with $\Gamma$ infinite property (T) group and $K$ arbitrary. Also,
if $\Gamma\in w\Cal T_0$ then $(\Gamma*K_0) \times K_1\in w\Cal
T_0$, $\forall K_0$ arbitrary $K_1$ infinite. The groups $G=\Bbb
Z^2 \rtimes \Gamma$, with $\Gamma \subset SL(2,\Bbb Z)$ non
amenable are w-rigid (cf. [Ka], [Ma], [Bu]), and so are all the
arithmetic groups in ([Va], [Fe]). However, these groups are not
in the class $w\Cal T_0$ because in each of these cases the
relatively rigid subgroup ($\Bbb Z^N$, for some $N$) is abelian.

\vskip .1in \noindent {\it 7.0.5. Definition.} An action $\sigma$
of $G$ on $A\simeq L^\infty(X, \mu)$ is {\it malleable} (1.5 in
[Po4]) if there exists a continuous action $\alpha : \Bbb R
\rightarrow {\text{\rm Aut}}(A \overline{\otimes} A)$ that
commutes with the ``double'' action $\sigma_g \otimes \sigma_g,
g\in G,$ and ``flips'' at some point $A \otimes 1$ onto $1 \otimes
A$. It is {\it s-malleable} if there also exists $\beta \in
{\text{\rm Aut}}(A \overline{\otimes} A)$ leaving $A=A \otimes 1$
pointwise fixed and satisfying $\beta^2=id$, $\beta \alpha_t =
\alpha_{-t}\beta, \forall t$. Also, $\sigma$ is {\it sub
malleable} (resp. {\it sub s-malleable}) if it can be
appropriately extended to a malleable (resp. s-malleable) action.
For the precise definitions see (1.5 and 4.2 in [Po4]).

\vskip .05in

Given an arbitrary infinite discrete group $G$, an example of an
action $\sigma$ of $G$ that is both sub s-malleable and clustering
is the (left) Bernoulli shift action with ``base space''
$(Y_0,\nu_0)$ ($Y_0 \neq$ single point). Recall that such $\sigma$
acts on $(X,\mu)=\Pi_g (Y_0,\nu_0)_g$ by
$\sigma_h((x_g)_g)=(x'_g)_g$, where $x'_g=x_{h^{-1}g}$, $\forall
h,g\in G$. Recall also that a generalized Bernoulli shift action
$\sigma$ of $G$ is defined as follows: Let $S$ be an infinite
countable set and $\pi$ a faithful action of $G$ on $S$; let
$(X,\mu)=\Pi_s (Y_0,\nu_0)_s$ and define the action $\sigma$ of
$G$ on $(X,\mu)$ by $\sigma_h((x_s)_s)=(x'_s)_s$, where
$x'_s=x_{\pi(h^{-1})(s)}$, $\forall h \in G$, $s\in S$. All
generalized Bernoulli shift actions are sub s-malleable. They are
mixing iff for any finite subset $S_0 \subset S$ there exists
$F\subset G$ finite such that $\pi(g)S_0 \cap S_0 =\emptyset$,
$\forall g \in G\setminus F$. This condition is easily seen to
also imply $\sigma$ is clustering.

To state the results, let us also recall that if $\gamma$ is a
character of $G$ then it defines an automorphism $\theta^\gamma$
of $M =L^\infty(X,\mu) \rtimes_\sigma G$ by $\theta(\Sigma_g a_g
u_g)=\Sigma_g \gamma(g)a_g u_g$, which is outer whenever $\gamma$
is non-trivial (see e.g. [T]).

\proclaim{7.1. Theorem (vNE Strong Rigidity)} Let $\sigma : G
\rightarrow {\text{\rm Aut}} (X, \mu)$, $\sigma_0 : G_0
\rightarrow {\text{\rm Aut}} (X_0, \mu_0)$ be free, ergodic
actions of countable discrete groups $G, G_0$ on standard
probability spaces $(X,\mu)$, $(X_0,\mu)$. Denote
$M=L^\infty(X,\mu) \rtimes_\sigma G$, $M_0 = L^\infty(X_0,\mu_0)
\rtimes_{\sigma_0} G_0$. Assume:

\vskip .05in  $(a)$. $G$ is ICC; $\sigma$ is sub malleable
clustering, e.g. $\sigma$ is a Bernoulli $G$-action;

$(b)$. $G_0$ w-rigid ICC. More generally, $G_0$ has a w-rigid,
wq-normal subgroup $H\subset G_0$ such that $\{hgh^{-1}\mid h \in
G_0\}$ infinite $\forall e\neq g \in G_0$ and $\sigma_{0|H}$
ergodic. (Note that $G_0 \in w\Cal T_1$.)

\vskip .05in

If $\theta: M_0 \simeq pMp$ is an isomorphism, for some projection
$p\in M$, then $p=1$ and there exist a character $\gamma$ of $G$
and isomorphisms $\delta : G_0 \simeq G$, $\Delta : (X_0, \mu_0)
\simeq (X,\mu)$ satisfying $\sigma(\delta(h)) \Delta = \Delta
\sigma_0(h)$, $\forall h\in G_0$, such that $\theta = {\text{\rm
Ad}}(u) \circ \theta^\gamma \circ \theta^{\delta,\Delta}$ for some
unitary $u \in M$.
\endproclaim
\vskip .05in \noindent {\it Proof}. Note first that if $H \subset
G_0$ is as in condition $(b)$ then $L(H)'\cap M_0 =\Bbb C$. This
is because the condition $\{hgh^{-1}\mid h \in G_0\}$ infinite
$\forall e\neq g \in G_0$ implies $L(H)'\cap M_0 \subset
L^\infty(X_0, \mu_0)$ and since $H$ acts ergodically on $X_0$ this
further implies $L(H)'\cap M_0 \subset (L^\infty(X_0,
\mu_0))^H=\Bbb C$.

We may clearly assume $p\in L(G)$. By (5.3 in [Po4]), there exists
a unitary $u_1\in pMp$ such that $u_1\theta(L(G_0))u_1^*\subset
pL(G)p$. We may thus assume, for simplicity of notations, that
$M_0=pMp$ with $pL(G)p \supset L(G_0)$. This allows us to apply
Theorem 5.1. Since $G$ ICC implies $G$ has no non-trivial finite
normal subgroups, it follows by 5.1 that $p=1$ and there exists a
unitary element $u_2 \in L(G)$ such that Ad$(u_2)(u^0_h)=
\gamma(\delta(h)) u_{\delta(h)}$, for some character $\gamma$ of
$G$ and an isomorphism of groups $\delta:G_0 \simeq G$. But this
is equivalent to the required conclusions. \hfill Q.E.D.

\proclaim{7.1'. Theorem} With the same notations as in $7.1$, if
we assume:

\vskip .05in  $(a)$. $G$ is ICC; $\sigma$ is sub s-malleable
clustering, e.g. $\sigma$ is a Bernoulli $G$-action;

$(b)$. $G_0\in w\Cal T_0$.

\vskip .05in

\noindent then the same conclusion as in $7.1$ holds true for any
isomorphism $\theta: M_0 \simeq pMp$.
\endproclaim
\vskip .05in \noindent {\it Proof}. The proof is exactly the same
as for 7.1, but using ($5.3'$ in [Po4]) instead of ($5.3$ in
[Po4]). \hfill Q.E.D.

Taking in 7.1 (resp. 7.1') both groups $G_0, G$ in $w\Cal T_1$
(resp. in $w\Cal T_0$) and both actions $\sigma_0, \sigma$ to be
Bernoulli shifts, we get a 1 to 1 functor  $G \mapsto
L^\infty(X,\mu)\rtimes_\sigma G=M$ from the class $w\Cal T_0\cup
w\Cal T_1$ (with group isomorphisms) to II$_1$ factors (with
algebra isomorphisms). Moreover, the same is true for the
associated II$_\infty$ factors $M^\infty \overset \text{\rm def}
\to =M\overline{\otimes} \Cal B(\ell^2(\Bbb N))$.

\proclaim{7.2. Corollary} Let $G_i\in w\Cal T_0 \cup w\Cal T_1$
(e.g. w-rigid ICC) and $\sigma_i : G_i \rightarrow {\text{\rm
Aut}} (X_i, \mu_i)$ be Bernoulli shift actions, $i=0,1$. Then the
$\text{\rm II}_1$ factors $M_i=L^\infty(X_i,\mu_i)
\rtimes_{\sigma_i} G_i$ (resp. $\text{\rm II}_\infty$ factors
$M_i^\infty$), $i=0,1$, are isomorphic if and only if the groups
$G_i$ are isomorphic and the actions $\sigma_i$ are conjugate,
$i=0,1$. Even more so, if $\theta: M^\infty_0 \simeq M^\infty_1$,
then $\theta$ is the amplification of an isomorphism $\theta_0 :
M_0 \simeq M_1$ of the form $\theta_0 = \text{\rm Ad}(u) \circ
\theta^\gamma \circ \theta^{\delta,\Delta}$, for some $u \in \Cal
U(M_1)$, $\gamma\in \text{\rm Char}(G_1)$ and $\delta:G_0 \simeq
G_1$, $\Delta:(X_0,\mu_0) \simeq (X_1,\mu_1)$ satisfying
$\sigma_1(\delta(h)) \Delta = \Delta \sigma_0(h)$, $\forall h\in
G_0$.
\endproclaim
\vskip .05in

For a Bernoulli shift action $\sigma$ of an ICC group $G$,
$L^\infty(X,\mu) \rtimes_\sigma G$ can be viewed as the canonical
``group measure space factor''-version of the ``group factor''
$L(G)$. Thus, Corollary 7.2 above can be regarded as an
affirmative answer to a natural ``relative variant'' of Connes'
rigidity conjecture ([C2]), stating that any isomorphism between
factors arising from ICC property (T) groups ought to come from an
isomorphism of the groups. Note though that this ``relative
variant'' holds true in a generality which, of course, fails to be
true for group algebras $L(G)$. For instance given ANY ICC group
$K$ (even $K$ amenable !) the group $G=SL(3,\Bbb Z) \times K$ is
w-rigid, so the factors $L^\infty(X,\mu) \rtimes G$ are
non-isomorphic for non-isomorphic $K$, by 7.2, while the
corresponding group factors $L(G)$ are mutually isomorphic if $K$
are taken amenable.

To state the next consequence, we need to recall some simple
general facts about automorphism groups associated with actions,
equivalence relations and von Neumann algebras.

One denotes by Aut$(\Cal R_\sigma)$ the group of automorphisms of
$(X,\mu)$ that preserve the orbits of $\sigma$ and by Inn$(\Cal
R_\sigma)$ the subgroup of automorphism whose graph belongs to
$\Cal R_\sigma$. Equivalently, $\alpha \in {\text{\rm Inn}}(\Cal
R_\sigma)$ iff there exists $u$ in the normalizer of
$A=L^\infty(X,\mu)$ in $M$ such that $\alpha=\text{\rm
Ad}(u)_{|A}$. Note that Aut$(\Cal R_\sigma)$ is the normalizer of
Inn$(\Cal R_\sigma)$ in ${\text{\rm Aut}} (X, \mu)$. We let
Out$(\Cal R_\sigma)=\text{\rm Aut}(\Cal R_\sigma)/\text{\rm
Int}(\Cal R_\sigma)$,

We denote by $\text{\rm Aut}(\sigma,G)$ the group of automorphisms
$\Delta\in \text{\rm Aut}(X,\mu)$ for which there exists an
(necessarily unique) automorphism $\delta=\delta(\Delta)$ of $G$
such that $\Delta \sigma_g \Delta^{-1} = \sigma_{\delta(g)},
\forall g$. A particularly important subgroup of $\text{\rm
Aut}(\sigma,G)$ is the commutant of $\sigma(G)$ in Aut$(X,\mu)$,
i.e. the group of automorphisms $\Delta \in \text{\rm
Aut}(\sigma,G)$ for which $\delta(\Delta)=id_G$. We denote this
subgroup by $\text{\rm Aut}_0(\sigma,G)$. It is clearly normal in
Aut$(\sigma,G)$.

It is not hard to see that if $\Delta \in \text{\rm
Aut}(\sigma,G)$ then $\Delta$ belongs to Inn$(\Cal R_\sigma)$ if
and only if $\Delta=\sigma_g$ for some $g \in G$ (see e.g. [MoSo]
or [Fu3]). Thus, the image of $\text{\rm Aut}(\sigma,G)$ in
Out$(\Cal R_\sigma)$ coincides with the quotient $\text{\rm
Out}(\sigma,G)\overset \text{\rm def} \to = \text{\rm
Aut}(\sigma,G)/\sigma(G)$.

Note that if $G$ is ICC then Aut$_0(\sigma,G)$ is isomorphic to
its image in $\text{\rm Aut}(\sigma,G)/\sigma(G)$ and thus to its
image in Out$(\Cal R_\sigma)$ as well. Note also that if $\sigma$
is a Bernoulli shift then there is a canonical group embedding
$\text{\rm Aut}(G)\ni \delta \mapsto \Delta_\delta \in \text{\rm
Aut}(\sigma,G)$ given by $\Delta_\delta ((x_g)_g)= (x'_g)_g$,
where $x'_g = x_{\delta^{-1}(g)}, \forall g\in G$. Thus, one has
the split exact sequence
$$
1 \rightarrow \text{\rm Aut}_0(\sigma,G) \rightarrow \text{\rm
Aut}(\sigma,G)) \rightarrow  \text{\rm Aut}(\sigma,G))/\text{\rm
Aut}_0(\sigma,G)=\text{\rm Aut}(G) \rightarrow 1
$$
which in case $G$ is ICC gives the split exact sequence:
$$
1 \rightarrow \text{\rm Aut}_0(\sigma,G) \rightarrow \text{\rm
Aut}(\sigma,G))/\sigma(G) \rightarrow \text{\rm Out}(G)
\rightarrow 1
$$

Finally, note that if we put $M=L^\infty(X,\mu) \rtimes_\sigma G$
and denote as usual Out$(M)=\text{\rm Aut}(M)/\text{\rm Int}(M)$,
then we have a natural embedding $\text{\rm Aut}(\Cal R_\sigma)
\subset \text{\rm Aut}(M)$ whose image in $\text{\rm Out}(M)$ is
Out$(\Cal R_\sigma)$. Another remarkable subgroup of $\text{\rm
Aut}(M)$ related to the action $\sigma$ is the group
$\{\theta^\gamma \mid \gamma \in \text{\rm Char}(G)\}\simeq
{\text{\rm Char}}(G)$, with $\theta^\gamma \in \text{\rm Aut}(M)$
as defined at the beginning of this section. This group is
isomorphic to its image in Out$(M)$ and it is normalized by the
image of Aut$(\sigma,G)$ in Aut$(M)$, with their intersection
being trivial.

With these remarks in mind, and denoting as before $M^\infty =
M\overline{\otimes} \Cal B(\ell^2(\Bbb N))$ the type II$_\infty$
factor associated with $M$, Corollary 7.2 trivially implies:

\proclaim{7.3. Corollary} With the above notations, assume $G\in
w\Cal T_0 \cup w\Cal T_1$ (e.g. $G$ w-rigid ICC) and $\sigma:G
\rightarrow \text{\rm Aut}(X,\mu)$ is a Bernoulli $G$-action. Then
$\text{\rm Out}(M^\infty)=\text{\rm Out}(M)= \text{\rm
Aut}_0(\sigma,G) \times (\text{\rm Char}(G) \rtimes \text{\rm
Out}(G))$. In particular, $\mycal F(M)=\{1\}$. Also, $\text{\rm
Out}(\Cal R_\sigma) = \text{\rm Aut}_0(\sigma,G) \times \text{\rm
Out}(G)$ and $\mycal F(\Cal R_\sigma)=\{1\}$.
\endproclaim
\vskip .05in

The above corollary  reduces the calculation of the outer
automorphism group Out$(M)$ of II$_1$ factors $M=L^\infty(X,\mu)
\rtimes_\sigma G$ (resp. equivalence relations $\Cal R_\sigma$)
with $G$ ICC group  in the class $w\Cal T_0\cup w\Cal T_1$ (e.g.
for $G$ w-rigid ICC group) and $\sigma$ Bernoulli shift of base
$(Y_0,\mu_0)$, to the calculation of the commutant
Aut$_0(\sigma,G)=\{\theta \in \text{\rm Aut}(X,\mu) \mid \theta
\sigma_g = \sigma_g \theta, \forall g\in G\}$ of the left
Bernoulli shift $\sigma$. We conjecture that any $\theta$
commuting with $\sigma$ must be of the form $\theta=\theta_0
\theta_1$ with $\theta_0$ a right shift by some element of the
group $G$ and $\theta_1$ a product type action $\theta_1=\Pi_g
\alpha_g$, where $\alpha_g = \alpha$, $\forall g$, for some
$\alpha \in \text{\rm Aut}(Y_0,\mu_0)$. Note that if this would be
the case, then Aut$_0(\sigma,G) \simeq G \times \text{\rm
Aut}(Y_0,\mu_0)$, so by 7.3 it would follow that Out$(M)= G \times
\text{\rm Aut}(Y_0,\mu_0) \times (\text{\rm Char}(G) \rtimes
\text{\rm Out}(G))$ and Out$(\Cal R_\sigma)=G \times \text{\rm
Aut}(Y_0,\mu_0) \times \text{\rm Out}(G)$. Such a calculation can
be shown to imply that two Bernoulli shift actions $\sigma_0,
\sigma_1$ of ICC groups $G_0, G_1\in w\Cal T_0\cup w\Cal T_1$ with
base space $(Y_0,\nu_0)$, $(Y_1,\nu_1)$ give rise to isomorphic
II$_1$ group measure space factors $M_0, M_1$ iff $G_0 \simeq G_1$
and $(Y_0,\nu_0) \simeq (Y_1,\nu_1)$.

The ``weak-normality'' assumptions on the relatively rigid
subgroup of $G_0$ in 7.1.(b), 7.1'.(b) allowed us to obtain plain
conjugacy of actions from von Neumann isomorphism (vNE). Using the
full generality of the Cartan Conjugacy Criteria in Section 4, one
can still obtain significant information (notably OE equivalence
from vNE equivalence) even if we drop such weak normality
assumption:

\proclaim{7.4. Theorem (vNE/OE Strong Rigidity)} Let $\sigma : G
\rightarrow {\text{\rm Aut}} (X, \mu)$, $\sigma_0 : G_0
\rightarrow {\text{\rm Aut}} (X_0, \mu_0)$ be free ergodic m.p.
actions. Denote $A=L^\infty(X,\mu), A_0=L^\infty(X_0, \mu_0)$, $M=
A \rtimes_\sigma G$, $M_0 = A_0 \rtimes_{\sigma_0} G_0$. Assume:

\vskip .05in  $(a)$. $G$ is ICC; $\sigma$ is sub s-malleable
clustering, e.g. $\sigma$ Bernoulli $G$-action. \vskip .05in
$(b)$. $G_0$ has a non virtually abelian subgroup $H \subset G_0$
with the relative property $(\text{\rm T})$.

\vskip .05in

If $\theta: M_0 \simeq pMp$ is an isomorphism, for some projection
$p\in M$, then $\tau(p)^{-1}$ is an integer and there exists a
partial isometry $v \in pMp$ such that

$(i)$. $v^*v = p$, $vv^* \in A$, $v\theta(A_0)v^* = Avv^*$.

$(ii)$. If $G_0'\subset G_0$ is the wq-normalizer of $H$ in $G_0$,
then there exists $G'\subset G$ and $K \vartriangleleft G'$,
$|K|=\tau(p)^{-1}$, such that $\text{\rm Ad}(v) \circ \theta$
conjugates $({\sigma_0})_{|G_0'}$ onto the action $\sigma'$ of
$G'/K$ on $A^K=\{ a\in A \mid \sigma_k(a)=a, \forall k \in K\}$
given by $\sigma'({\hat{g'}}) = (\sigma_{g'})_{|A^K}$, $\forall
g'\in G'$.

\endproclaim
\vskip .05in \noindent {\it Proof}. Like in the proofs of 7.1,
7.1', by ([Po4]) we may assume $p \in L(G)$ and $\theta(L(G'_0))
\subset pL(G)p$. By Theorem 4.2 it follows that $\theta(A_0)$ is
unitary conjugate to $Ap_0$ for some projection $p_0 \in A$ with
$\tau(p_0)=\tau(p)$. But then 5.2 applies to get the conclusion.
\hfill Q.E.D.

The first result deriving  OE (equivalently, Cartan conjugacy)
from vNE, i.e. from the isomorphism of the group measure space
factors, was obtained in ([Po2]): It is shown there that if
$M_0=M(\Cal R_0)$, $M=\Cal M(\Cal R)$ for some equivalence
relations $\Cal R_0, \Cal R$, with $\Cal R_0$ having the relative
property (T) (as defined in 5.10.1 in [Po2]) and $\Cal R$ having
the Haagerup property (as considered in 2.1 and 3.5.6 of [Po2])
then any isomorphism $\theta:M_0 \simeq M$ is an inner
perturbation of an isomorphism coming from OE. Once such a
``vNE/OE Strong Rigidity'' is proved, one can take advantage of
results on OE Ergodic Theory to derive results about the II$_1$
factors involved. Thus, in ([Po2]) one uses Gaboriau's
$\ell^2$-Betti numbers for equivalence relations to derive
invariants for the corresponding group measure space factors,
(called the {\it class} $\Cal H\Cal T$).

Similarly, one can use 7.4 in combination with ([Ga2]), or with
the Monod-Shalom OE Strong Rigidity results in ([MoSh]), to obtain
calculation of fundamental groups and non-isomorphism of II$_1$
group measure space factors that are not covered by Theorem 7.1.
For instance, one can apply 7.4 to the case $G_0, G$ are products
of two or more groups of the form $H*K$ with $H,K$ torsion free
and $H$ Kazhdan, to get OE from vNE, then apply ([MoSh]) to
further get conjugacy from OE, and the triviality of the
fundamental group. Along these lines, we mention a consequence of
7.4 which derives the triviality of fundamental group of certain
factors without using OE rigidity results from ([Ga2], [MoSh]):

\proclaim{7.5. Corollary} If $G$ is an ICC group having a non
virtually abelian subgroup with the relative property $(\text{\rm
T})$ (e.g. $G=H*K$ for some infinite Kazhdan group $H$ and $K\neq
1$, or a product of such groups) and $\sigma$ is a Bernoulli
$G$-action, then the fundamental group of the corresponding group
measure space factor $M$ satisfies $\mycal F(M)\cap [1, \infty)
\subset \Bbb N$. If in addition $G$ is torsion free (e.g. $G=H *
K$ with $H$ Kazhdan and $H,K$ torsion free) then $\mycal
F(M)=\{1\}$.
\endproclaim
\vskip .05in

Since any OE of actions entails an isomorphism of the associated
group measure space factors, Theorems 7.1, 7.1' trivially imply OE
rigidity results for actions $\sigma, \sigma_0$ satisfying
conditions $(a), (b)$ of either 7.1 or 7.1'. But in fact we can do
slightly better than that: Since an isomorphism of group measure
space factors that comes from an OE isomorphism of probability
spaces takes Cartan subalgebras onto each other, we can use 5.2
(which only requires mixing conditions on the actions) rather than
5.1 (which is what we used in the proofs of 7.1, 7.1', and which
required the clustering condition on $\sigma$). This allows us to
replace the condition ``$\sigma$ clustering'' by the weaker
condition ``$\sigma$ mixing''.

\proclaim{7.6. Theorem (OE Strong Rigidity)} Let $\sigma : G
\rightarrow {\text{\rm Aut}} (X, \mu)$, $\sigma_0 : G_0
\rightarrow {\text{\rm Aut}} (X_0, \mu_0)$ be free ergodic m.p.
actions. Assume:

\vskip .05in  $(a)$. $G$ is ICC; $\sigma$ is sub s-malleable and
mixing (e.g. a Bernoulli $G$-action). \vskip .05in $(b)$. Either
$G_0 \in w\Cal T_0$; or $(\sigma_0, G_0)$ satisfies $7.1.(b)$
(e.g., $G_0$ w-rigid ICC)

\vskip .05in

If $Y\subset X$ is a subset of positive measure and $\Delta_0:
(X_0, \mu_0) \simeq (Y,\mu_Y)$ satisfies $\Delta_0(\Cal
R_{\sigma_0})=\Cal R^Y_\sigma$ then $Y=X$, the groups $G_0,G$ are
isomorphic and the actions $\sigma_0, \sigma$ are conjugate with
respect to that isomorphism of groups. More precisely, there exist
$\alpha \in \text{\rm Inn}(\Cal R_\sigma)$, $\delta : G_0 \simeq
G$, $\Delta : (X_0, \mu_0) \simeq (X,\mu)$ satisfying
$\sigma(\delta(h)) (\Delta(x)) = \Delta(\sigma_0(h)(x))$, $\forall
h\in G_0, x\in X_0$, such that $\Delta_0 = \alpha \circ \Delta$.
\endproclaim
\vskip .05in \noindent {\it Proof}. With the usual notation $M_0,
M$ for the associated group measure space factors, the isomorphism
$\Delta_0$ implements an isomorphism $\theta : M_0 \simeq pMp$
satisfying $\theta(A_0)=Ap$, where $p=\chi_Y$. Since $L(G)$ is
diffuse (because $G$ is infinite), $p$ can be unitary conjugate to
a projection in $L(G)$. Equivalently, after identifying $M_0$ with
$\theta(M_0)$, we may assume $M_0=pMp$ for some projection $p\in
L(G)$. By applying ($5.3$ in [Po4]), in the case $(\sigma_0, G_0)$
satisfies $7.1.(b)$, and respectively ($5.3'$ in [Po4]) in the
case $G_0 \in w\Cal T_1$, it follows that we may also assume
$L(G_0) \subset pL(G)p$. Thus, assumptions $(a)$ and $(c)$ of 5.2
are satisfied. Since $G$ is ICC and $M_0=pMp$, all the conditions
required in the last part of 5.2 are checked and the corresponding
conclusion in 5.2 trivially implies the statement. \hfill Q.E.D.

\vskip .05in Following ([MoSh]), to have an ``OE
superrigidity''-type result one needs all conditions on just one
side. Theorem 7.6, which has the conditions on the group (either
w-rigid ICC or in the class $w\Cal T_0$) on the ``source'' side
and the conditions on the action (Bernoulli) on the ``target''
side, doesn't quite fulfill  this requirement. By strengthening
the condition on the group though, we can obtain this kind of
result as well:

\proclaim{7.7. Theorem (OE Superrigidity)} Let $G$ be an ICC
property $(\text{\rm T})$ group and $\sigma: G \rightarrow
\text{\rm Aut}(X,\mu)$ a Bernoulli $G$-action (more generally, a
sub s-malleable, mixing action). Let $G_0$ be any group and
$\sigma_0 : G_0 \rightarrow \text{\rm Aut}(X_0, \mu_0)$ any free
m.p. action of $G_0$. If $\Delta_0: (X_0, \mu_0) \simeq
(Y,\mu_{|Y})$ satisfies $\Delta_0(\Cal R_{\sigma_0})=\Cal
R^Y_\sigma$, for some subset $Y\subset X$, then $Y=X$ $($modulo
null sets$)$, the groups $G_0,G$ are isomorphic and the actions
$\sigma_0, \sigma$ are conjugate with respect to that isomorphism
of groups. More precisely, there exist $\alpha \in \text{\rm
Inn}(\Cal R_\sigma)$, $\delta : G_0 \simeq G$ such that $\Delta =
\alpha \circ \Delta_0$ satisfies $\sigma(\delta(h)) \Delta =
\Delta \sigma_0(h)$, $\forall h\in G_0$.
\endproclaim
\vskip .05in \noindent {\it Proof}. By either (4.1.7, 4.1.9 in
[Po8]) or (1.4 in [Fu2]), if $\Cal R_{\sigma_0} \simeq \Cal
R^Y_\sigma$ and $G$ has property (T) then $G_0$ has property (T).
But then we can apply 7.5 and the statement follows. \hfill Q.E.D.

\vskip .05in Note that the class of groups  and actions covered by
the OE rigidity results 7.6 and 7.7 is quite different from the
ones in ([Fu1], [MoSh]). However, some intersection occurs with
([MoSh]). For instance, if we take $G_0=\Gamma_1 \times ... \times
\Gamma_n$ with $n \geq 2$, each $\Gamma_i$ either a free product
of torsion free groups or a non-elementary Gromov-hyperbolic group
with the property (T), at least one being of this latter class,
then each $\Gamma_i$ is in the class $\Cal C$ of Monod-Shalom
while $G_0$ is weakly rigid. Thus, if the action $\sigma_0$ of
$\Gamma_1 \times ... \times \Gamma_n$ is assumed irreducible, in
the sense of ([MoSh]), and $\sigma$ is Bernoulli with $G$ ICC,
then the corresponding OE rigidity result is covered by both 7.5
above and (1.9 in [MoSh]). On the other hand, we were not able to
check Furman's condition for Bernoulli shift actions $\sigma$ of
higher rank lattices $\Gamma\subset \Cal G$, i.e. to show whether
$(\sigma, G; X)$ has no quotients of the form $\Cal G/\Lambda$,
with $\Lambda$ a lattice in $\Cal G$. If this would be true (which
is what we expect), then (Theorem A in [Fu1]) would imply: Let $G$
be an ICC higher rank lattice and $\sigma$ Bernoulli $G$-action.
If $\Cal R_\sigma^Y \simeq \Cal R_{\sigma_0}$ for some $Y \subset
X$ and some  free ergodic m.p. action $(\sigma_0, G_0)$, then $G,
G_0$ are ``virtually isomorphic'' groups and $\sigma, \sigma_0$
``virtually conjugate'' actions.

Our next result gives a Galois type correspondence for the w-rigid
ICC sub-equivalence relations of an equivalence relation
implemented by a Bernoulli shift. The result can also be viewed as
a OE Strong Rigidity result for embeddings of equivalence
relations.

\proclaim{7.7. Theorem (OE Strong Rigidity for Embeddings} Let
$\sigma : G \rightarrow {\text{\rm Aut}} (X, \mu)$ be a Bernoulli
$G$-action (or more generally a sub s-malleable mixing m.p.
action) of an ICC group $G$. If $\sigma_0:G_0 \rightarrow
{\text{\rm Aut}} (X_0, \mu_0)$ is a free ergodic m.p. action of  a
group $G_0$ which is either w-rigid ICC or in the class $w\Cal
T_0$ and $\Delta_0 : (X_0, \mu_0) \simeq (X,\mu)$ satisfies
$\Delta_0(\Cal R_{\sigma_0})\subset \Cal R_\sigma$ then there
exists an isomorphism $\delta: G_0 \simeq G'\subset G$ and $\alpha
\in {\text{\rm Inn}}(\Cal R_\sigma)$ such that $\alpha \circ
\Delta_0 \circ \sigma_0(h) = \sigma(\delta(h))$, $\forall h\in
G_0$. Thus, there is a $1$ to $1$ correspondence between the
sub-equivalence relations of $\Cal R_\sigma$ that are
implementable by free ergodic actions of w-rigid ICC groups (resp.
groups in $w\Cal T_0$), modulo conjugacy by automorphisms in
$\text{\rm Inn}(\Cal R_\sigma)$, and the w-rigid ICC subgroups of
$G$ (resp. the $w\Cal T_0$ subgroups of $G$), modulo conjugacy by
elements in $G$.
\endproclaim
\vskip .05in \noindent {\it Proof}. We may assume $\Cal
R_{\sigma_0}\subset \Cal R_\sigma$, i.e. we can select unitaries
$\{u^0_h\}_{h\in G_0}$ in the normalizer of $A=L^\infty(X,\mu)$ in
$M=L^\infty(X,\mu) \rtimes_\sigma G$ such that $G_0 \ni h \mapsto
\sigma_0(h)=\text{\rm Ad}(u^0_h) \in \text{\rm Aut}(X,\mu)$,
$\forall h\in G_0$. Since $A \subset M=M(\Cal R_\sigma)$ has
vanishing 2-cohomology in the sense of ([FM]), it follows that for
any sub-equivalence $\Cal R \subset \Cal R_\sigma$ the Cartan
inclusion $A \subset M(\Cal R)$ has vanishing 2-cocycle. Thus, by
perturbing if necessary each $u^0_h$ by a unitary element in $\Cal
U(A)$ we may assume $\{u^0_h\}_{h\in G_0}$ are the canonical
unitaries generating $L(G_0)$. If $G_0$ is w-rigid ICC then by
(5.3 in [Po4]) there exists a unitary element $u$ in $M$ such that
$uL(G_0)u^* \subset L(G)$. Since $\sigma$ is clustering, by Lemma
4.5 it follows that $\sigma_0$ is mixing, thus Theorem 5.2 can be
applied to conclude that there exists $u\in \Cal U(M)$ such that
$u(\{u^0_h\}_h)u^* \subset \Bbb T \{u_g\}_g$ and $uAu^*=A$. The
proof of the case $G_0 \in w\Cal T_0$ is similar, using ($5.3'$ in
[Po4]). \hfill Q.E.D.

\vskip .05in In (part $(1)$ of Theorem D in [Fu1]), Alex Furman
proved that the restrictions to subsets $Y$ of irrational measure
of the equivalence relation $\Cal R_\sigma$ given by a classical
action $\sigma$ of an arithmetic Kazhdan group, such as $SL(n,\Bbb
Z)$ on $\Bbb T^n$, $n \geq 3$, cannot be implemented by a free
action of a group. We end by noticing that Theorem 7.7 provides a
new class of examples of this situation, without restrictions on
the measure of the subsets. Note however that (Thm. D$(3)$ in
[Fu1]) also provides examples of equivalence relations $\Cal R$
with the property that no amplification $\Cal R^t, t > 0$ can be
implemented by a free action of a group, a situation which we
cannot realize here (see however Remark 5.6.2 in [Po4]) !

\proclaim{7.8. Corollary} Let $\sigma : G \rightarrow {\text{\rm
Aut}} (X, \mu)$ be a Bernoulli shift action (or more generally a
free, sub s-malleable, mixing  m.p. action) of an  ICC property
$(\text{\rm T})$ group $G$. If $Y \subset X$ has measure $0\neq
\mu(Y) < 1$ then $\Cal R^Y_\sigma$ (the restriction of the
equivalence relation $\Cal R_\sigma$ to $Y$) cannot be realized as
orbits of a free action of a group.
\endproclaim

\head References\endhead

\item{[Bu]} M. Burger: {\it Kazhdan constants for} SL$(3,\Bbb Z)$,
J. reine angew. Math., {\bf 413} (1991), 36-67.

\item{[Chr]}
E. Christensen: {\it Subalgebras of a finite algebra},
Math. Ann. {\bf 243} (1979), 17-29.

\item{[C1]} A. Connes: {\it A type II$_1$
factor with countable fundamental group}, J. Operator Theory {\bf
4} (1980), 151-153.

\item{[C2]} A. Connes: {\it Classification des facteurs}, Proc.
Symp. Pure Math. {\bf 38} (Amer. Math. Soc. 1982), 43-109.

\item{[C3]} A. Connes: {\it Classification of injective factors},
Ann. of Math., {\bf104} (1976), 73-115.

\item{[C4]} A. Connes: {\it Sur la classification des facteurs de
type} II, C. R. Acad. Sci. Paris {\bf 281} (1975), 13-15.

\item{[CJ1]} A. Connes, V.F.R. Jones: {\it A} II$_1$ {\it factor
with two non-conjugate Cartan subalgebras}, Bull. Amer. Math. Soc.
{\bf 6} (1982), 211-212.

\item{[CJ2]} A. Connes, V.F.R. Jones: {\it Property} (T) {\it for
von Neumann algebras}, Bull. London Math. Soc. {\bf 17} (1985),
57-62.

\item{[CSh]} A. Connes, D. Shlyakhtenko: {\it $L^2$-homology for
von Neumann algebras}, preprint,
math.OA/0309343.

\item{[CoHa]} M. Cowling, U. Haagerup: {\it Completely bounded
multipliers of of the Fourier algebra of a simple Lie group of
real rank one}, Invent. Math. {\bf 96} (1989), 507-549.

\item{[CoZi]} M. Cowling, R. Zimmer: {\it Actions of lattices in}
$Sp(n,1)$, Ergod. Th. Dynam. Sys. {\bf 9} (1989), 221-237.

\item{[Dy]} H. Dye: {\it On groups of measure preserving
transformations}, II, Amer. J. Math, {\bf 85} (1963), 551-576.

\item{[FM]} J. Feldman, C.C. Moore: {\it Ergodic equivalence relations,
cohomology, and von Neumann algebras I, II}, Trans. Amer. Math.
Soc. {\bf 234} (1977), 289-324, 325-359.

\item{[Fe]} T. Fernos: {\it Kazhdan's Relative Property} (T): {\it
Some New Examples}, math.GR/0411527.

\item{[Fu1]} A. Furman: {\it Orbit equivalence rigidity}, Ann. of
Math. {\bf 150} (1999), 1083-1108.

\item{[Fu2]} A. Furman: {\it Gromov's measure equivalence and
rigidity of higher rank lattices}, Ann. of Math. {\bf 150} (1999),
1059-1081.

\item{[Fu3]} A. Furman: {\it Outer automorphism groups of some
ergodic equivalence relations}, preprint 2003.

\item{[Ga1]} D. Gaboriau: {\it Cout des r\'elations d'\'equivalence
et des groupes}, Invent. Math. {\bf 139} (2000), 41-98.

\item{[Ga2]} D. Gaboriau: {\it Invariants $\ell^2$ de r\'elations
d'\'equivalence et de groupes},  Publ. Math. I.H.\'E.S. {\bf 95}
(2002), 93-150.

\item{[GeGo]} S. Gefter, V.  Golodets: {\it Fundamental groups for
ergodic actions and actions with unit fundamental groups}, Publ.
RIMS, Kyoto Univ. {\bf 24} (1988), 821-847.

\item{[GoNe]} V. Y. Golodets, N. I. Nesonov: T{\it -property and
nonisomorphic factors of type} II {\it and} III, J. Funct.
Analysis {\bf 70} (1987), 80-89. \item{[Gr]} M. Gromov: {\it
Hyperbolic groups} in ``Essays in group theory'' (S.M.S. Gersten
editor), Springer, 1987, 75-263.

\item{[Ha]} U. Haagerup: {\it An example of a non-nuclear
C$^*$-algebra which has the metric approximation property},
Invent. Math. {\bf 50} (1979), 279-293.

\item{[Hj]} G. Hjorth: {\it A converse to Dye's theorem},
UCLA preprint, Sept. 2002.

\item{[dHV]} P. de la Harpe, A. Valette: ``La propri\'et\'e T
de Kazhdan pour les groupes localement compacts'', Ast\'erisque
{\bf 175} (1989).

\item{[J1]} V.F.R. Jones : {\it Index for subfactors}, Invent.
Math. {\bf 72} (1983), 1-25.

\item{[J2]} V.F.R. Jones : {\it Ten problems}, in ``Mathematics:
perspectives and frontieres'', pp. 79-91, AMS 2000, V. Arnold, M.
Atiyah, P. Lax, B. Mazur Editors.

\item{[Ka]} D. Kazhdan: {\it Connection of the dual space of a group
with the structure of its closed subgroups}, Funct. Anal. and its Appl.,
{\bf1} (1967), 63-65.

\item{[Ma]} G. Margulis: {\it Finitely-additive invariant measures
on Euclidian spaces}, Ergodic. Th. and Dynam. Sys. {\bf 2} (1982),
383-396.

\item{[MoSh]} N. Monod, Y. Shalom:
{\it Orbit equivalence rigidity and bounded cohomology},
Preprint 2002.

\item{[Moo]} C.C. Moore: {\it Ergodic theory and von Neumann
algebras}, Proc. Symp. Pure Math. {\bf 38} (Amer. Math. Soc.
1982), 179-226.

\item{[MvN1]} F. Murray, J. von Neumann:
{\it On rings of operators}, Ann. Math. {\bf 37}
(1936), 116-229.

\item{[MvN2]} F. Murray, J. von Neumann:
{\it Rings of operators IV}, Ann. Math. {\bf 44}
(1943), 716-808.

\item{[NeSt]} S. Neshveyev, E. St\o rmer: {\it Ergodic
theory and maximal abelian subalgebras of the hyperfinite
factor}, JFA {\bf 195} (2002), 239-261.

\item{[Oz]} N. Ozawa: {\it A Kurosh type theorem for type} II$_1$
{\it factors}, math.OA/0401121.

\item{[Po1]} S. Popa: {\it Some rigidity results for non-commutative
Bernoulli shifts}, MSRI preprint 2001-005.

\item{[Po2]} S. Popa: {\it On a class of
type II$_1$ factors with Betti numbers invariants}, MSRI preprint
2001-024, revised math.OA/0209310, to appear in Ann. of Math.

\item{[Po3]} S. Popa: {\it On the fundametal group
of type} II$_1$ {\it factors}, Proc. Nat. Acad. Sci. {\bf 101}
(2004), 723-726. (math.OA/0210467).

\item{[Po4]} S. Popa: {\it Strong rigidity of} II$_1$ {\it factors
arising from malleable actions of weakly rigid groups}, I,
math.OA/0305306.

\item{[Po5]} S. Popa: {\it On the distance between MASA's in type}
II$_1$ {\it factors}, Fields Institute Comm. {\bf 30} (2001),
321-324.

\item{[Po6]} S. Popa: {\it Some computations of $1$-cohomology groups
and construction of non orbit equivalent actions}, math.OA/0407199.

\item{[Po7]} S. Popa: {\it Orthogonal pairs of $*$-subalgebras in
finite von Neumann algebras}, Journal of Operator Theory, {\bf 9}
(1983), 253-268.

\item{[Po8]} S. Popa: {\it Correspondences}, INCREST preprint
1986.

\item{[Po9]} S. Popa: {\it Classification of subfactors of type}
II, Acta Math. {\bf 172} (1994), 163-255.

\item{[PoSa]} S. Popa, R. Sasyk: {\it On the cohomology of actions of
groups by Bernoulli shifts}, math.OA/0311417.

\item{[PoSiSm]} S. Popa, A. Sinclair, R. Smith: {\it Perturbations of
subalgebras of type} II$_1$ {\it factors}, J. Funct. Analysis {\bf
213} (2004), 346-379 (math.OA/0305444).

\item{[S]} I.M. Singer: {\it Automorphisms of finite factors},
Amer. J. Math. {\bf 177} (1955), 117-133.

\item{[T]} M. Takesaki: ``Theory of Operator Algebras II'',
Encyclopedia of Mathematical Sciences {\bf 125}, Springer-Verlag,
Berlin-Heidelberg-New York, 2002.

\item{[Va]} A. Valette: {\it Group pairs with relative property} (T)
{\it from arithmetic lattices}, preprint 2004 (preliminary version
2001).

\item{[Zi]} R. Zimmer: ``Ergodic theory and semisimple groups'',
Birkha\"user-Verlag, Boston 1984.

\enddocument